\documentclass[reqno]{amsart}
\usepackage{fancyhdr}
\usepackage{manfnt}
\usepackage{qtree}
\usepackage{amsmath,amsthm,amscd,amsfonts,amssymb,euscript}
\usepackage{graphics}
\usepackage{pdfpages}
\usepackage{mathtools}

\usepackage[all]{xypic}
\usepackage{latexsym}
\usepackage{color} 
\swapnumbers
\theoremstyle{plain}

\DeclareFontFamily{OT1}{rsfs}{}
\DeclareFontShape{OT1}{rsfs}{n}{it}{<-> rsfs10}{}
\DeclareMathAlphabet{\mathscr}{OT1}{rsfs}{n}{it}

\newcommand{\eE}{{\mathscr E}}
\newcommand{\eG}{{\mathscr G}}
\newcommand{\eH}{{\mathscr H}}

\begin{document}
\input{amssym.def}


\newtheorem{theorem}[subsubsection]{\sc Theorem}
\newcommand{\bth}{\begin{theorem}$\!\!\!${\bf }~}
\newcommand{\eeth}{\end{theorem}}

\newtheorem{propo}[subsubsection]{\sc Proposition}

\newcommand{\bprop}{\begin{propo}$\!\!\!${\bf }~}
\newcommand{\eprop}{\end{propo}}

\newtheorem{lema}[subsubsection]{\sc Lemma}
\newcommand{\blem}{\begin{lema}$\!\!\!${\bf }~}
\newcommand{\elem}{\end{lema}}

\newtheorem{defe}[subsubsection]{\sc Definition}
\newcommand{\bdefe}{\begin{defe}$\!\!\!${\bf}~}
\newcommand{\edefe}{\end{defe}}

\newtheorem{coro}[subsubsection]{\sc Corollary}
\newcommand{\bcor}{\begin{coro}$\!\!\!${\bf }~}
\newcommand{\ecor}{\end{coro}}

\theoremstyle{remark}
\newtheorem{rema}[subsubsection]{\it Remark}
\newcommand{\brem}{\begin{rema}$\!\!\!${\it }~\rm}
\newcommand{\erem}{\end{rema}}

\newtheorem{note}[subsubsection]{Note}

\newtheorem{observe}[subsubsection]{\it Observation}
\font\manual=manfnt
\def\dbend{{\manual\char127}} 

\def\danger{\begin{trivlist}\begin{footnotesize}\item[]\noindent%
\begingroup\hangindent=3pc\hangafter=-2
\def\par{\endgraf\endgroup}%
\hbox to0pt{\hskip-\hangindent\dbend\hfill}\ignorespaces}
\def\enddanger{\par\end{footnotesize}\end{trivlist}}

\def\ddanger{\begin{trivlist}\begin{footnotesize}\item[]\noindent%
\begingroup\hangindent=3pc\hangafter=-2
\def\par{\endgraf\endgroup}%
\hbox to0pt{\hskip-\hangindent\dbend\kern2pt\dbend\hfill}\ignorespaces}
\def\endddanger{\par\end{footnotesize}\end{trivlist}}

\newcommand{\spec}{{\rm Spec}\,}
\newtheorem{notation}{Notation}[subsubsection]
\newcommand{\bnot}{\begin{notation}$\!\!\!${\bf }~~\rm}
\newcommand{\enot}{\end{notation}}
\numberwithin{equation}{subsubsection} 
\newcommand{\beqa}{\begin{eqnarray}}
\newcommand{\eeqa}{\end{eqnarray}}
\newtheorem{thm}{Theorem}
\theoremstyle{definition}
\newtheorem{hint}[thm]{Hint}
\newtheorem{example}[subsubsection]{\it Example}
\newtheorem{say}[subsubsection]{}
\newcommand{\bsem}{\begin{say}$\!\!\!${\it }~~\rm}
\newcommand{\esem}{\end{say}}

\newsymbol \circledarrowleft 1309

\newcommand{\ha}{\sf h}
\newcommand{\g}{\sf g}
\newcommand{\ta}{\sf t}
\newcommand{\s}{\sf s}
\newcommand{\ctext}[1]{\makebox(0,0){#1}}
\setlength{\unitlength}{0.1mm}

\newcommand{\wt}{\widetilde}
\newcommand{\Lr}{\Longrightarrow}
\newcommand{\Aut}{\mbox{{\rm Aut}$\,$}}
\newcommand{\ul}{\underline}
\newcommand{\ol}{\bar}
\newcommand{\lr}{\longrightarrow}
\newcommand{\bc}{{\mathbb C}}
\newcommand{\bp}{{\mathbb P}}
\newcommand{\bz}{{\mathbb Z}}
\newcommand{\bq}{{\mathbb Q}}
\newcommand{\bn}{{\mathbb N}}
\newcommand{\bg}{{\mathbb G}}
\newcommand{\br}{{\mathbb R}}
\newcommand{{\bh}}{{\mathbb H}}

\newcommand{\cl}{{\mathcal L}}
\newcommand{\cv}{{\mathcal V}}
\newcommand{\cf}{{\mathcal F}}
\newcommand{\cb}{{\Lambda}}
\newcommand{\mfc}{{\sf C}}
\newcommand{\ce}{{\mathcal E}}
\newcommand{\co}{R}
\newcommand{\cs}{\mathcal S}
\newcommand{\cg}{{\mathcal G}}
\newcommand{\ca}{{\mathcal A}}
\newcommand{\hra}{\hookrightarrow}
\newcommand{\mfu}{{\sf U}}

\newtheorem{ack}{\it Acknowledgments}       
\renewcommand{\theack}{} 

\title[Parahoric bundles]{Moduli of parahoric $\mathcal G$--torsors on a compact Riemann surface} 
\author{V. Balaji} 
\author{C. S. Seshadri} 
\address{Chennai Mathematical Institute SIPCOT IT Park, Siruseri-603103, India,
balaji@cmi.ac.in, css@cmi.ac.in}
\date{}
\dedicatory{Dedicated to Professor K. Chandrasekharan in admiration.}
\thanks{The research of the first author was partially supported by the J.C. Bose Research grant.}
\begin{abstract}
Let $X$ be an irreducible smooth projective algebraic curve of genus $g \geq 2$ over the ground field $\bc$ and let $G$ be a semisimple simply connected algebraic group. The aim of this paper is to introduce the notion of  {\em semistable and stable parahoric}  torsors under a certain Bruhat-Tits group scheme $\mathcal G$ and construct the moduli space of semistable parahoric $\mathcal G$--torsors; we also identify the underlying topological space of this moduli space with certain spaces of homomorphisms of Fuchsian groups into  a maximal compact subgroup of $G$.  The results give a generalization of the earlier results of Mehta and Seshadri on parabolic vector bundles. 
\end{abstract}
\maketitle

\vspace{2mm}
\noindent
\section{Introduction}  Let $X$ be a  smooth projective curve defined over $\bc$ of genus $g \geq 2$. Let ${\mathcal R} \subset X$ be a fixed set of points of $X$  with $m = |{\mathcal R}|$ and let $n_i$ be a set of positive integers attached to each of the points $x_i \in {\mathcal R}$. The  uniformization theorem states that there exists a simply connected covering  surface $q:\tilde X \to X$, unique upto isomorphism, subordinate to the given  signature i.e. ramified precisely over the points ${\mathcal R} = \{x_i\} \subset X$ together with ramification indices $n_i$ at these points.  Since $g \geq 2$, we may identify $\tilde X$ with $\bh$ the {\em upper half space} (cf. \cite[page 49-50]{weil}). Let $\pi$ be the subgroup of the discontinuous group of automorphisms of $\bh$  such that $X = {\bh}/{\pi}$. Note that the action of $\pi$ is \textit{not} free. Let $q:\bh \to X$ be the quotient projection.   It is well known that the isotropy subgroups at the points $z_i \in q^{-1}({\mathcal R})$ are cyclic of finite order. Let  the isotropy subgroups be denoted by 
\[
\pi_{z_i} = \langle C_i \rangle 
\]
with $C_i$ as generators.
 Thus each $C_i$ is an element of order $n_i$.

We fix once for all the set ${\mathcal R}$ and the indices $n_i$. Recall that the group $\pi$  is a Fuchsian group generated by $2g + m$ elements $A_1, B_1,\ldots, A_g, B_g, C_1,\ldots, C_m$, modulo the relations
\beqa\label{intro5}
A_{_1}B_{_1}A_{_1}^{-1}B_{_1}^{-1} \ldots A_{_g}B_{_g}A_{_g}^{-1}B_{_g}^{-1}\ldots C_{1}
\ldots C_{m} \, = I.
\eeqa
\beqa\label{intro5.1} C_{_i}^{n_i} = I, ~~~~(i=1,2, \ldots,m).
\eeqa

Let $G$ be a connected reductive algebraic group over $\bc$ and let $K_G \subset G$ be a maximal compact subgroup of $G$. 
\bdefe\label{localtypeofrep} The {\em type of a homomorphism} $\rho:\pi \to G$ is defined to be the set of conjugacy classes in $G$ of the images $\rho(C_i)$ and is denoted by ${\boldsymbol\tau}= \{\boldsymbol\tau_i \}$. Equivalently, the type of $\rho$ is the set of isomorphism classes of the local representations $\rho_{{z_i}}:\pi_{z_i} \to G, i = 1, \ldots, m$.    \edefe 

\bnot Let $R^{\boldsymbol\tau}(\pi, K_G)$ denote the space of homomorphisms $\rho:\pi \to K_G$ of type $\boldsymbol\tau = \{\boldsymbol\tau_i\}$. \enot

\bdefe\label{pigee} A  $(\pi,G)$--bundle on $\bh$  is defined to be the trivial $G$--bundle $\bh \times G$ on $\bh$ with the  $\pi$--structure given by $\gamma(z,g) = (z, \rho(\gamma).g)$, with $\rho$ a homomorphism $\pi \to G$.\edefe 

If $G = GL(n)$ is the full-linear group, the $(\pi,G)$--bundles on $\bh$ have an equivalent description as  $\pi$--vector bundles on $\bh$. We recall (\cite{pibundles}, \cite{ms}) that if $V \simeq \bh \times {\bc}^n $ is a $\pi$--vector bundle on $\bh$, the vector bundle $W = q^{\pi}_{_*}(V)$(invariant direct image by $q$) on $X$ acquires a {\em parabolic structure} which consists of the data assigning a {\em flag to the fibre of $W$ at every ramification point in $X$ for the covering $q$} together with a tuple of {\em weights}. 

The invariant direct image functor $V\mapsto q_*^\pi(V)$ gives a fully faithful
embedding of the category of $\pi$--vector bundles on ${\mathbb H}$ into the
category of parabolic vector bundles on $X$ (morphisms being taken as isomorphisms).
Moreover, we can realise every parabolic bundle with rational weights as $q_*^\pi(V)$ for a suitable $\pi$ and $V$ (cf. \cite{ms}).

This translates easily into an equivalent description of  $(\pi,GL(n))$--bundles on $\bh$ as principal $GL(n)$--bundles on $X$ with {\em parabolic structures}. One can define the concepts of {\em stability} (resp. {\em semistability}) for $\pi$--vector bundles (or equivalently parabolic bundles on $X$) and construct the corresponding moduli space of equivalence classes of semistable objects (fixing some invariants) as a normal projective variety.  As topological spaces  these moduli spaces can be identified with set of equivalence classes of elements in $R^{\boldsymbol\tau}(\pi, U(n))$, i.e. {\em unitary representations of $\pi$} (see Mehta-Seshadri \cite{ms}, Seshadri \cite{pibundles}), which generalize the results in Narasimhan-Seshadri\cite{ns} and Seshadri \cite{unitary}. 

The purpose of this paper is to generalize the above results when the structure group $G$ is no longer the full-linear group. Let us suppose hereafter that the group $G$ is {\em semisimple and simply connected} (over $\bc$) unless otherwise stated.

 One can again give an equivalent description of $(\pi, G)$--bundles on $\bh$ as certain intrinsically defined objects on $X$. However, the picture is more subtle than the case when $G$ is the full-linear group ; for instance, {\em it is not possible, in general,  to associate in a natural manner a principal $G$--bundle on $X$ to a $(\pi,G)$--bundle on $\bh$}. The new objects on $X$, which give an equivalent description of $(\pi,G)$--bundles on $\bh$, will be called parahoric bundles or {\em parahoric torsors}. These parahoric torsors are defined as pairs $(\mathscr E, \boldsymbol\theta)$, where $\mathscr E$  is a torsor (i.e. principal homogeneous space) on $X$ under a parahoric Bruhat-Tits group scheme $\mathcal G$, together with a prescription of {\em weights} $\boldsymbol\theta$, which are elements of the set of {\em rational} one-parameter subgroups of $G$ (see the discussion below and Definition \ref{quasiparahoric}). We define notions of semistability and stability of such parahoric torsors and  construct  moduli spaces of these objects. 
   
 The torsors under parahoric group schemes that we consider here have been studied earlier by Pappas and Rapoport, without however the notion of weights (see \cite{pr1} and \cite{pr2}); in \cite{pr2} they made some precise conjectures on the moduli stack of such torsors. Heinloth has since settled many of their conjectures (see \cite{heinloth}; we note that Heinloth works over arbitrary ground fields not just $\bc$). We were led to the study of  parahoric torsors in trying to interpret $(\pi,G)$--bundles on $\bh$ as objects on $X$ (inspired by A. Weil's work \cite{weil}, as was the case in \cite{ms} and \cite{pibundles}). In Section 2 we  link explicitly the ideas from the paper of Weil and Bruhat-Tits theory. This relationship plays a key role in the rest of the paper. We need to define a few technical terms before we can state the main results of our paper.

Let $A_{{x_i}}$ be the completion of the local ring at $x_i$,  and $K_{{x_i}}$ (or simply as $K$) its quotient field , $x_i \in {\mathcal R}$. Let $T$ be a maximal torus of $G$ and $Y(T) := Hom({\mathbb G}_m, T)$ the group of $1$--parameter subgroups of $T$; let ${\mathbb E} \simeq Y(T) \otimes {\br}$ and ${\mathbb E}_{\bq} \simeq Y(T) \otimes \bq$.  By the general theory of Bruhat and Tits (see \cite[Definition 5.2.6]{bt}) and \ref{boun} below), one has  certain collection of subsets $\{\Theta_i\} \subset {\mathbb E}_{\bq}^m$, where $m = |{\mathcal R}|$ and to each subset $\Theta_i \subset {\mathbb E}$, one can associate a parahoric subgroup ${\mathcal P}_{_{{\Theta}_i}}(K) \subset G(K_{{x_i}})$, $i = 1, \ldots m$, and furthermore, associated to each parahoric subgroup ${\mathcal P}_{_{{\Theta}_i}}(K)$, there is a smooth group scheme ${\mathcal G}_{_{\Theta_i}}$ over $D_{x_i} = Spec~A_{{x_i}}$, known as a Bruhat-Tits group scheme.

More precisely, (see \ref{boun}), by fixing a root datum the theory of buildings allows us to identify the vector space ${\mathbb E}$ with an {\em affine apartment} $App(G,K)$ in the Bruhat-Tits building, and each parahoric subgroup ${\mathcal P}_{_{{\Omega}_i}}(K) \subset G(K)$ is precisely the stabilizer subgroup of a {\em facet} $\Omega_i$ of the affine apartment for the natural $G(K)$--action on the building. We will reserve the symbol $\Omega$ for a facet of the apartment.

 It is explained in Section 5 as to how, given any finite subset of points ${\mathcal R} \subset X$ and a collection ${\mathcal G}_{_{\Theta_i}}$ of {\em Bruhat-Tits group schemes} over $D_{x_i}$, one can construct a (global) group scheme ${\mathcal G}_{_{\Theta,X}}$ over the projective curve $X$ by {\em gluing} (see Lemma \ref{gettingbt} and Definition \ref{globalbt}) so that
\beqa{\mathcal G}_{_{\Theta,X}}|_{_{_{X - {\mathcal R}}}} \simeq G \times (X - {\mathcal R}), ~~~~~~ {\mathcal G}_{_{\Theta,X}}|_{D_{x_i}} \simeq {\mathcal G}_{_{\Theta_i}}, x_i \in {\mathcal R}. \eeqa
We will call the set ${\mathcal R}$ the points of ramifications of  $\mathcal G$.  
 
 Following Pappas and Rapoport (\cite{pr2}), we will call ${\mathcal G}_{_{\Theta,X}}$ the {\bf parahoric Bruhat-Tits group schemes}.  However, we wish to emphasize that both Pappas-Rapoport (\cite{pr2}) and Heinloth \cite{heinloth}, {\em do not} make the assumption that  ${\mathcal G}_{_{\Theta,X}}|_{_{_{X - {\mathcal R}}}} \simeq G \times (X - {\mathcal R})$, i.e., for them the group scheme need not be {\em generically  split}.

 It can be shown (see Remark \ref{assocbtgrpscheme}) that to every set $\boldsymbol\tau$ of conjugacy classes and finite subset $\mathcal R \subset X$, we can associate a collection ${\boldsymbol\theta_{\boldsymbol\tau}} = \{\theta_i\} \in {\mathbb E}_{\bq}^m $  of elements of ${\mathbb E}_{\bq}$ and also a parahoric  group scheme ${\mathcal G}_{_{\boldsymbol\theta_{{\boldsymbol\tau},X}}}$ on $X$ such that the points of ramifications of ${\mathcal G}_{_{\boldsymbol\theta_{{\boldsymbol\tau},X}}}$ is $\mathcal R$. The content of Theorem \ref{maintheoremI} below is that this correspondence $\boldsymbol\tau \mapsto {\mathcal G}_{_{\boldsymbol\theta_{{\boldsymbol\tau},X}}}$ extends precisely to give an identification of moduli spaces of representations with fixed conjugacy classes and that of torsors under ${\mathcal G}_{_{\boldsymbol\theta_{{\boldsymbol\tau},X}}}$.

One of the key features of parahoric groups is that for any interior point $\theta$ of a facet $\Omega_i$, we have an isomorphism ${\mathcal P}_{_{{\Omega_i}}}(K) \simeq {\mathcal P}_{_{\theta}}(K)$ (see the discussion in \ref{boun} below). In particular, any parahoric Bruhat-Tits group scheme ${\mathcal G}_{_{\Omega,X}}$ associated to a collection of facets $\{\Omega_i\}$ is isomorphic  to a ${\mathcal G}_{_{\boldsymbol\theta_{{\boldsymbol\tau},X}}}$ for some $\boldsymbol\tau$.

Before going to the main results of this paper, we begin by observing that a collection ${\boldsymbol\theta_{\boldsymbol\tau}} = \{\theta_i\}$ of rational weights   entails a choice of ramification indices $d_i$ at the points of $\mathcal R$ (see Remark \ref{ramindexandtheta}). Since the genus $g \geq 2$, it is well-known (see \ref{ramindices}),  that there exists a Galois cover $p:Y \to X$, with Galois group $\Gamma$, ramified precisely at $\mathcal R$ with the prescribed ramification indices $d_i$. 

It is shown in Theorem \ref{stackygammaversusparahoric} that there is an  isomorphism between the moduli stack of $(\Gamma,G)$--bundles on $Y$ of local type $\boldsymbol\tau$ and the stack  of ${\mathcal G}_{_{\boldsymbol\theta_{{\boldsymbol\tau},X}}}$--torsors on $X$.


We then define, in Section 6 of this paper, the concept of {\em semistable and stable} $\mathcal G$--torsors on $X$ as well as the notion of $S$--equivalence. 
Our main results can be formulated as follows ( see Theorem \ref{maintheoremI,II}, Theorem \ref{realdimension} and Corollary \ref{dimensionofmodulispace} for notation and details):
\bth\label{maintheoremI}
Let ${\mathcal G}_{_{\boldsymbol\theta_{{\boldsymbol\tau},X}}}$ be a parahoric Bruhat-Tits group scheme associated to $\boldsymbol\tau$. 
\begin{enumerate}{\it 
\item The set $M_{_X}({\mathcal G}_{_{\boldsymbol\theta_{{\boldsymbol\tau},X}}})$ of $S$--equivalence classes of semistable ${\mathcal G}_{_{\boldsymbol\theta_{{\boldsymbol\tau},X}}}$--torsors on $X$  gets a natural structure of an irreducible normal projective variety of dimension 
\beqa
dim_{\bc}(G)(g -1) + \sum_{i = 1}^{m}{\frac{1}{2}}e({\boldsymbol\theta_{\boldsymbol\tau}})
\eeqa
In fact, the variety $M_{_X}({\mathcal G}_{_{\boldsymbol\theta_{{\boldsymbol\tau},X}}})$ is the coarse moduli space for the functor of isomorphism classes of ${\mathcal G}_{_{\boldsymbol\theta_{{\boldsymbol\tau},X}}}$--torsors on $X$. 
\item Let ${\overline K_G} = K_G/centre$. There exists a Fuchsian group $\pi$  and a bijective correspondence between the space $R^{\boldsymbol\tau}(\pi, K_G)/{{\overline K_G}}$ of conjugacy classes of homomorphisms $\rho:\pi \to K_G$ of local type $\boldsymbol\tau$ and the set of $S$--equivalence classes of {\em semistable} ${\mathcal G}_{_{\boldsymbol\theta_{{\boldsymbol\tau},X}}}$--torsors. 
\item This correspondence induces a homeomorphism  
\[R^{\boldsymbol\tau}(\pi, K_G)/{{\overline K_G}}\simeq M_{_X}({\mathcal G}_{_{\boldsymbol\theta_{{\boldsymbol\tau},X}}})
\] 
of the underlying topological spaces.
\item Under this correspondence, the subset of irreducible homomorphisms gets identified with isomorphism classes of  {\em stable} ${\mathcal G}_{_{\boldsymbol\theta_{{\boldsymbol\tau},X}}}$--torsors.}

\end{enumerate}
\eeth

We make a few clarifying remarks on the paper.

\brem
\noindent
\begin{enumerate}

\item The moduli stack ${\sf Bun}_{_X}({\mathcal G}_{_{\boldsymbol\theta_{{\boldsymbol\tau},X}}})$ has been studied in detail by Heinloth (\cite{heinloth}).

\item ({\em Parabolic $G$--bundles})~   If for a point $x \in {\mathcal R}$ the parahoric group ${\mathcal P}_{_{{\Omega}}}(K_x)$ gets identified with the distinguished hyperspecial parahoric subgroup $G(A_x)$ (see \ref{hyperspecial} for the definition)  the moduli space of parahoric torsors gets identified with the moduli space of principal $G$--bundles on $X$ in the usual sense. In this case, the parahoric structure comes from the origin of $\mathbb E$ (see  \ref{boun}).

If ${\mathcal P}_{_{{\Omega}}}(K_x)$ is a proper subgroup of  $G(A_x)$  then under the evaluation map $ev: G(A_x) \to G(\bc)$, the subgroup ${\mathcal P}_{_{{\Omega}}}(K_x)$ is the inverse image of a standard parabolic subgroup of $G$, so that in this case a quasi-parahoric torsor (see Definition \ref{parahorictorsor}) could indeed be called a {\em quasi-parabolic $G$--bundle} in the familiar sense of the term when $G = GL(n)$ is the full-linear group, i.e. the data consists of a principal $G$--bundle on $X$ together with a parabolic subgroup of $G$ (i.e. a ``flag") for every $x \in {\mathcal R}$. This case corresponds to the situation when the parahoric subgroups defining the local Bruhat-Tits group schemes come from the interior of the Weyl alcove. Equivalently the {\em weights} come from the interior of the Weyl alcove. These are the cases dealt with in Teleman-Woodward (\cite{tw}).

\item ({\em Parahoric torsors which are not principal $G$--bundles}) We now consider parahoric subgroups of $G(K_x)$ which cannot be conjugated to subgroups of $G(A_x)$. For instance, barring $G(A_x)$, the rest of the maximal parahoric subgroups of $G(K_x)$ fall under this case (see \cite{bt}). The weights in these cases lie on the walls of the Weyl alcove (cf. Teleman \cite[Section 9]{teleman}).

 It is this case which highlights one of the reasons why we need to give a subtler description of $(\Gamma,G)$--bundles on $Y$ as parahoric torsors on $X$ which do not support a principal $G$--bundle on $X$. Evidence to this effect was shown using Tannakian considerations in Balaji-Biswas-Nagaraj \cite{bbn1}, leading to the definition of a ramified bundle in \cite{bbn2}.  More concrete examples were shown in \cite{sesramanan} indicating what to expect in general.

\item  The striking cases which arise out of the present study are the non-hyperspecial maximal parahoric subgroups where a number of new phenomena show up. These correspond, on the side of the representations of the Fuchsian group (see \ref{intro5}), to  those maps $\rho:\pi \to K_G$ such that centralizers of the images of the elements $\rho(C_i)$ are  {\em proper semisimple} subgroups of $G$ (see Remark \ref{assocbtgrpscheme}).

\item After this paper was posted in the archives, we were informed by P. Boalch of his paper \cite{boalch} where the parahoric structure is seen in the setting of regular singular connections. 
\end{enumerate}
\erem

{\it Acknowledgements}: We wish to thank Jochen Heinloth and Michel Brion for many helpful suggestions on an earlier version of this paper. The first author also thanks Gopal Prasad  for some  helpful discussions on Bruhat-Tits theory. The first author wishes to thank the Isaac Newton Institute for their hospitality during the semester on ``Moduli" where this work was given the final shape. We wish to thank Pramathanath Sastry and Brian Conrad for their very helpful comments. Finally we sincerely thank the  referee for the meticulous reading of the paper and  comments and questions which have helped immensely in clarifying many key issues in the paper.

\section{Non-abelian functions and bounded groups}

\subsection~As the title suggests, the aim of this section is to tie up some ideas from the classical paper of A. Weil (\cite{weil}) and Seshadri (\cite{pibundles}) and Bruhat-Tits theory (\cite{bt}). This section is central to this paper. 
\subsubsection{ \it Some preliminaries on root data}
Let $G$ be a semisimple, simply connected algebraic group  defined over $\bc$; we fix a maximal torus 
$T$ of $G$. Let $X(T) := Hom(T, {\mathbb G}_m)$ be the character group and $Y(T) := Hom({\mathbb G}_m, T)$ the group of $1$--parameter subgroups of $T$. Let  $R = R(T,G) \subset {X(T)}$ be the root system  associated to the adjoint representation of $G$ and  $S$ be a  system of simple roots. 

Denote by $(~,~):Y(T) \times X(T) \to {\bz}$ the canonical bilinear form. The set $S$ determines a system of positive roots $R^{+} \subset {R}$ and a Borel subgroup $B \subset G$ with unipotent radical $U$. We now order the set $R^+ = \{r_i\}, i = 1, \ldots, q$. We then have a family $\{u_{_r}:{\bg}_a \to G \mid r \in R \}$ of {\em root homomorphisms} of groups such that one gets an isomorphism of varieties:
\beqa
\prod_{r \in R^+}u_{_{r}}:\prod_{r \in R^+} {\bg}_a \to U
\eeqa

 For every root $r \in R$, we denote by $T_{_r} = Ker(r)^0$, and $Z_{_r} = Z_G(T_{_r})$, the centralizer of $T_r$ in $G$. The derived group $[Z_{_r}, Z_{_r}]$ is of rank 1 and there exists  a unique 1PS, ${r^{\vee}}:{\bg}_m \to T \cap [Z_{_r}, Z_{_r}]$ such that $T = Im({r}^{\vee}).T_{_r}$ and $({r}^{\vee}, r) = 2$.
The element $r^{\vee}$ is the coroot (or 1--PS) associated to $r$. 

For each $r \in R$ the  {\em root homomorphism}
\beqa\label{root1}
u_{_r}:{\mathbb G}_a \to G
\eeqa
is such that
\beqa\label{root2}
t.u_{_r}(a).t^{-1} = u_{_r}(r(t).a)
\eeqa
for any $\bc$--algebra $A$, $t \in T(A)$, $a \in A$, and such that the tangent map $du_{_r}$ induces an isomorphism
\[
du_{_r}:Lie({\mathbb G}_a) \to (Lie G)_{_r}
\]
The functor $A \mapsto u_{_r}({\mathbb G}_a) = u_{_r}(A)$ gives $U_{_r}(A) \subset G(A)$. This determines a closed subgroup $U_{_r}$ of $G$ and is called the {\em root group} corresponding to $r$. 

Denote by $\{{\alpha}^* \mid \alpha \in S\}$ the coroots  dual to $\{\alpha \in S\}$, i.e. $(\alpha^*, r) = \delta_{\alpha,r}$. Define 
\beqa\label{e and e'}
{\mathbb E} := Y(T) \otimes_{\bz} {\br}
\\
{\mathbb E}':= X(T) \otimes_{\bz} {\br} 
\eeqa
Most often, we in fact work with $X(T) \otimes_{\bz} {\bq}$ and $Y(T) \otimes_{\bz} {\bq}$.

\subsubsection{\it Parahoric subgroups}\label{boun}  Let $K$ be  the field $\bc((z))$ of Laurent power series in $z$ and let $A = \bc[[z]]$ be the ring of integers, with residue field $\bc$. 

For the notion of Bruhat-Tits buildings and their behaviour under field extensions see J.Tits \cite[Page 43]{tits}.

Once we fix a root datum for $G$, we see that we have a choice of an affine apartment; the choice of the maximal torus $T$ then identifies ${\mathbb E}$ with an affine apartment $App(G,K)$ in the Bruhat-Tits building ${\mathcal B}(G,K)$.

A subset $M \subset G(K)$ is said to be {\em bounded} if for any
regular function $f \in K[G]$, the values $v(f(m))$ for the valuation $v$ on $A$, are bounded
below when $m$ runs over all elements of $M$. In particular, we may talk of {\em bounded subgroups}. A subgroup $M \subset G(K)$ is therefore bounded if the ``order of poles" of elements of $M$ is bounded. This can be made precise by taking a faithful representation  $G \hra GL(n)$ so that elements of $M$ are represented by matrices with entries in $K$. 

Let $\Theta \subset {\mathbb E}$ be a nonempty subset which is a facet.  Denote by 
${\mathcal P}_{_\Theta}(K) \subset G(K)$ the subgroup generated by $T(A)$ and the groups $U_{r}(z^{m_{r}} A)$ for all the roots $r \in R$, where 
\beqa\label{genparahorics}
m_{r} = m_{r}(\Theta) = -[inf_{_{\theta \in \Theta}} (\theta, r)]
\eeqa
where $[h]$ stands for the biggest integer smaller than or equal to $h$.

The group ${\mathcal P}_{_\Theta}(K)$  is a bounded subgroup, more precisely it is a {\em parahoric subgroup} of $G(K)$ in the sense of Bruhat-Tits and conversely, any parahoric subgroup is bounded in the above sense (cf. Bruhat-Tits \cite{bt}) .

The choice of a root datum  identifies a parahoric subgroup ${{\mathcal P}_{_\Omega}(K)} \subset G(K)$ as the stabilizer subgroup of $G(K)$ of a {\em facet} $\Omega$ of the affine apartment $App(G,K)$ for the natural $G(K)$--action on ${\mathcal B}(G,K)$. By Tits \cite[Section 3.1, page 50]{tits}, since we work with a semisimple and simply connected group $G$ we could in turn take any point in {\em general position} i.e. an interior point in the facet and consider the parahoric subgroup as the stabilizer of that point. Thus one can make an identification ${{\mathcal P}_{_\Omega}(K)} \simeq {{\mathcal P}_{_\theta}(K)}$ for an {\em interior point} $\theta$ in the facet  $\Omega$. 

By the main theorem of Bruhat-Tits (\cite{bt}), there exist smooth group schemes ${\mathcal G}_{_\Omega}$ over $Spec(A)$ such that the group ${\mathcal G}_{_\Omega}(A) =  {\mathcal P}_{_\Omega}(K)$ and moreover, since $A$ is a complete discrete valuation ring, the group scheme is uniquely determined upto unique isomorphism by its $A$--valued points (see \cite[Section 1.7]{bt}).

Let $\theta \in {\mathbb E}$. Thus,
\beqa
m_{r} = m_{r}(\theta) = -[(\theta, r)] 
\eeqa
In other words, we have:
\beqa\label{gille-1}
{\mathcal P}_{_\theta}(K) = \langle T(A),~~~~ U_{r}(z^{m_r(\theta)}A), ~~~r \in R \rangle.
\eeqa

To summarize, since we work with a semisimple and simply connected group $G$, all {\em parahoric groups} are, upto conjugacy by elements of $G(K)$, precisely the collection of groups $\{{\mathcal P}_{_\theta}(K)\}_{\theta \in {\mathbb E}}$ (see \cite[Section 3.1,~page 50]{tits}), and as such we will work with these groups. In fact, we may choose these $\theta$ to be in ${\mathbb E}_{\bq} = Y(T) \otimes \bq$. Again  by \cite[page 51]{tits}, the conjugacy classes of maximal parahoric subgroups of $G(K)$ are the stabilizers of the vertices of the building and they are precisely $l+1$ in number, where $l = rank(G)$. In particular, associated to the ``origin" $0 \in {\mathbb E}$ we have the group ${\mathcal P}_{_0}(K)$, which is nothing but the maximal bounded subgroup $G(A) \subset G(K)$. 

Note that if $\theta$ lies in the lattice $Y(T)$ itself, then there exists $t \in T(K)$ such that 
\beqa\label{gille0}
{\mathcal P}_{_\theta}(K) = t .{\mathcal P}_{_0}(K). t^{-1}
\eeqa
\brem\label{interiorofalcove1} Again we note that if $m_r(\theta) < 1$ for all $r \in R$, then ${\mathcal P}_{_\theta}(K) \subset G(A)$. These parahoric subgroups then correspond to the standard parabolic subgroups of $G$.\erem



\brem In this remark we make some comments on the parahoric groups when we make the assumption that the group $G$ is {\em simple}. Let $\alpha_{_{max}}$ denote the highest root. Then we can express it as:
\beqa\label{calpha}
\alpha_{_{max}} = \sum_{\alpha \in S} c_{\alpha} \cdot \alpha
\eeqa 
with $c_{\alpha} \in {\bz}^+$.

One can have a nicer choice of the points whose stabilizers give the maximal parahorics (see the last paragraph in \cite[Page 662]{tits1}), now that $G$ is {\em simple}. For every $\alpha \in S$, if we define 
\beqa\label{thetaalpha}
\theta_{\alpha} = {\frac {\alpha^*}{c_{\alpha}}} \in {\mathbb E},
\eeqa
{\em then in fact, the groups $\{{\mathcal P}_{_{\theta_{\alpha}}}(K) \}_{_{\alpha \in S}} \cup {\mathcal P}_{_0}(K) $  represent the conjugacy classes under $G(K)$ of all {\em maximal parahoric subgroups} of $G(K)$}. In other words, these are indexed precisely by the vertices of the {\em extended Dynkin diagram}.\erem

\subsubsection{\it Hyperspecial Parahorics} \label{hyperspecial}  In  Bruhat-Tits theory, we encounter the so-called {\em hyperspecial} maximal parahorics which have the following characterizing property: each parahoric group ${\mathcal P}_{_\Omega}(K)$ is identified with ${\mathcal G}_{_\Omega}(A)$, the $A$--valued points  of a certain canonically defined smooth group scheme ${\mathcal G}_{_\Omega}$ defined over $A$. It is a fact that the parahoric subgroup ${\mathcal P}_{_{\theta_{\alpha}}}(K)$ is hyperspecial if and only if $c_{\alpha} = 1$ in the description of the long root $\alpha_{_{max}}$. The hyperspecial parahorics are listed at the end of \cite{tits}.

\bsem{\em The Weyl alcove} We now recall the description of the set of conjugacy classes in a compact semisimple and simply connected group in terms of the affine Weyl group $W_{_{\text{aff}}}$. 

Let $K_G \subset G$ be a maximal compact subgroup. For an arbitrary group $S$, let $\text{Torsion}(S)$ denote the  subset of elements of finite order in $S$. We then have the following identifications:
\[
\text{Torsion}(K_G)/conjugation \simeq \text{Torsion}(T)/W 
\]
\[
(Y(T) \otimes {\bq}/{\bz})/W \simeq (Y(T) \otimes \bq) /W_{_{\text{aff}}}.
\] 
If further, the group $G$ is assumed to be {\em simple}, then $(Y(T) \otimes \bq) /W_{_{\text{aff}}}$ gets identified with the simplex ({\em the (rational) Weyl alcove}) 
\[
{\mathcal A}:= \{x \in Y(T) \otimes \bq \mid (x,\alpha_{_{max}}) \leq 1, (x, \alpha_i) \geq 0, \forall~~positive~~ roots~~\alpha_i  \}
\] 

\brem\label{alcove}  In fact,  the set of conjugacy classes of element in $K_G$ gets identified with $T/W$ which is the Weyl alcove since any element of $K_G$ is conjugate to an element in the maximal torus upto an element of the Weyl group (cf. \cite[page 151]{morgan}).    
\erem \esem

\brem\label{assocbtgrpscheme} Recall that vertices of the alcove ${\mathcal A}$ correspond to the vertices of the extended Dynkin diagram. Furthermore, to each point of ${\mathcal A}$ one can associate a {\em parahoric subgroup} of $G(K)$ and hence a canonically defined parahoric Bruhat-Tits group scheme. Thus, for each tuple ${\boldsymbol\tau} = \{{\boldsymbol\tau}_i \}_{i = 1}^{m}$ of conjugacy classes of elements of finite order in $K_G$ we have a point $\boldsymbol\theta_{\boldsymbol\tau} = \{{\theta}_i\}_{i = 1}^{m} \in {\mathcal A}^{m}$, where $m = \#\{ ~of~~ conjugacy~~ classes\}$ and hence an associated parahoric Bruhat-Tits group scheme ${\mathcal G}_{_{\boldsymbol\theta_{\boldsymbol\tau}}}$.

More can be said. Let $\alpha \in S$ is a simple root, and ${\mathcal G}_{_{\theta_{\alpha}}}$ the Bruhat-Tits group scheme associated to the maximal parahoric ${\mathcal P}_{_{\theta_{\alpha}}}(K)$. Let $g_\alpha$ be an element in $K_G$ of finite order corresponding to ${\theta_{\alpha}}$. Then the centralizer $Z_{_G}(g_{\alpha})$ can be obtained from the {\em closed fibre of the group scheme} ${\mathcal G}_{_{\theta_{\alpha}}}$; indeed, $Z_{_G}(g_{\alpha}) \simeq \big({\mathcal G}_{_{\theta_{\alpha}}}\big){_{_{_x}}}/{\{unipotent~radical\}}$. When ${\theta_{_\alpha}}$ is {\em hyperspecial}, then ${\mathcal G}_{_{\theta_{\alpha}}}$ is in fact a semisimple group scheme and therefore  
$Z_{_G}(g_{\alpha}) = G$. On the other hand, when ${\theta_{_\alpha}}$ is non-hyperspecial, $Z_{_G}(g_{\alpha})$'s are precisely those subgroups of $G$ which are proper semisimple subgroups of maximal rank in $G$ corresponding to the classical Borel-de Siebenthal list (see \cite{borel}). 

\erem

\brem\label{alcoveparahoric} In the case when $G$ is assumed to be simple and simply connected, by the description of the (rational) Weyl alcove $\mathcal A$ (see Definition \ref{alcove}) and the fact that the parahoric subgroups are determined by {\em interior points} of $\mathbb E$, it follows that upto conjugacy by $G(K)$, every parahoric subgroup of $G(K)$ can be identified with a ${\mathcal P}_{_\theta}(K)$ for a suitable $\theta \in {\mathcal A}$. Moreover, by Remark \ref{interiorofalcove1}, if $m_r(\theta) < 1$ for all $r \in R$, then ${\mathcal P}_{_\theta}(K) \subset G(A)$. \erem
\brem\label{semisimplealcove} We remark that even when $G$ is  semisimple, we still have the notion of an {\em alcove} ${\mathcal A}$, but it will no longer be a simplex as in the case when $G$ is simple since there is no unique $\alpha_{_{max}}$ but ${\mathcal A}$ will now be a product of the Weyl alcoves associated to the simple factors of $G$. Again, parahoric subgroups will be parametrized by points of the alcove upto conjugacy by $G(K)$. \erem

\subsubsection{\it Standard parahorics} \label{stdparahorics} (See Remark \ref{interiorofalcove1}, Remark \ref{interiorofalcove} and Remark \ref{alcoveparahoric}) Following the loop group terminology, the {\em standard parahoric subgroups} of $G(K)$ are parahoric subgroups of the distinguished hyperspecial parahoric subgroup  $G(A)$. These are realized as inverse images under the evaluation map $ev: G(A) \to G(k)$ of standard parabolic subgroups of $G$. For any $I \subset S$, let $R_{_I}$ denote the set $R_{_I} = R \cap {\mathbb Z}I$. Let $U_{_I} := U((-R^{+})\setminus R_I)$ and $L_{_I} :=G(R_I)$. The standard parabolic $P_{_I} \subset G$ is defined by $P_{_I}:= U_{_I} L_{_I}$.

In particular, the standard {\em Iwahori subgroup} ${\mathfrak I}$ is a standard parahoric  and indeed, ${\mathfrak I} = ev^{-1}(B)$, with $B = U(-R^{+}) T$ being the standard Borel subgroup  containing the fixed maximal torus $T$.

Since the {\em standard parahoric subgroups} of $G(A) = {\mathcal P}_{_{0}}(K)$ are also indexed by the subsets of the set of simple roots, to avoid any confusion, we will denote the {\em standard parahoric subgroups} of $G(A)$ by ${\mathcal P}^{st}_{_{I}}(K):= ev^{-1}(P_{_I})$ for every subset $I \subset S$. 

For instance if $\alpha \in S$ let $S_{_\alpha} := S \setminus \{\alpha\}$. Then $P_{_{S_{_\alpha}}} \subset G$ is a maximal parabolic subgroup and $ev^{-1}(P_{_{S_{_\alpha}}}) = {\mathcal P}^{st}_{_{S_{_\alpha}}}(K)$ is a standard parahoric which can be described as follows:
\beqa
{\mathcal P}^{st}_{_{S_{_\alpha}}}(K) = \langle T(A),~~~~ U_{r}(A), ~~~r \in R_{_{S_{_\alpha}}} \cup (-R^{+})\setminus R_{_{S_{_\alpha}}}\rangle
\eeqa
Note that
$ev^{-1}(L_{_{S_{_\alpha}}}) = \langle T(A),~~~~ U_{r}(A), ~~~r \in R_{_{S_{_\alpha}}} \rangle$ and 
$ev^{-1}(U_{_{S_{_\alpha}}}) = \langle U_r(A), r \in  (-R^{+})\setminus R_{_{S_{_\alpha}}}\rangle$.

If $r \in R_{_{S_{_\alpha}}}$,  the simple root $\alpha$ does not occur in $r$, in which case $(\theta_{\alpha},r) = 0$. Hence $ev^{-1}(L_{_{S_{_\alpha}}})  \subset {\mathcal P}_{_{\theta_{\alpha}}}(K) \cap {\mathcal P}_{_{0}}(K)$.

\noindent
Again if $r \in  (-R^{+})\setminus R_{_{S_{_\alpha}}}$, say $r = \sum a_{\beta}.\beta$, with $a_{\beta} \leq 0$ and $a_{\alpha} \neq 0$. By the definition of $c_{\alpha}$, we have $-1 \leq \frac{a_{\alpha}}{c_{\alpha}}  < 0$. It follows that $m_r(\theta_{\alpha}) = -[(\theta_{\alpha},r)] = -[\frac{a_{\alpha}}{c_{\alpha}}] = 1$. Hence~~$ev^{-1}(U_{_{S_{_\alpha}}}) \subset {\mathcal P}_{_{\theta_{\alpha}}}(K) \cap {\mathcal P}_{_{0}}(K)$.

We therefore have the inclusions:
\beqa\label{inclusions}
{\mathfrak I} \subset {\mathcal P}^{st}_{_{S_{_\alpha}}}(K) \subset {\mathcal P}_{_{\theta_{\alpha}}}(K) \cap {\mathcal P}_{_{0}}(K).
\eeqa

These standard parahorics will play a role when we re-look at Hecke correspondences.

\subsection{\bf Non-abelian functions and the unit group}

For the purposes of working in the  setting of algebraic curves instead of $\bh$, we make a few observations. A result due to A.Selberg (\cite{selberg}) states that if $A \subset GL(n,\bc)$ is a finitely generated subgroup, then $A$ has a {\em normal subgroup} $A_0$ of finite index with no torsion. It follows from this  that the discrete group $\pi \subset Aut(\bh)$ has a normal subgroup $\pi_0$ of finite index  such that $\pi_0$ operates \textit{freely} on ${\bh}$. Let $Y ={\bh} / \pi_0$ and $\Gamma = \pi/\pi_0.$ Then there is a canonical action of $\Gamma$ on $Y$ such that $X = Y/ \Gamma$. Let $p:Y \to X$ be the covering map and note that $\Gamma = Gal(Y/X)$. Conversely, if $Y$ is a Galois cover of $X$ with the  given signature (i.e., $\mathcal R$ and ramification indices), by the universality of $q:\bh \to X$ for this signature  it follows that there is a $\pi_o\subset \pi$ acting freely on $\bh$ such that $Y = \bh/{\pi_o}$.

\brem\label{ramindices}~~~~~In other words, given a finite number of points $x_i \in X$ together with signatures or ramification indices $n_i$ at these points, there exist ramified Galois covers $p:Y \to X$ albeit non-canonical, ramified precisely over the $x_i$ with the given ramification indices.  \erem

Let $r:{\bh} \to Y$ be the simply connected covering projection of $Y$. By definition, $\pi_o = \pi_1(Y)$ and we have the following commutative diagram:
\beqa\label{intro1}
\xymatrix{
{\bh} \ar[dr]_{q} \ar[rr]^{r}& & Y \ar[dl]^{p} \\
& X &
}
\eeqa
with $q = p \circ r$. Let $y_i$ be the image of $z_i$ in $Y$ and let ${\mathcal R} = \{x_i\}$, with $\{x_i = p(y_i) \mid  1 \leq i \leq m\}$.

 The map $r:{\bh} \to Y$ is a local isomorphism; in fact, if $z \in {\bh}$ maps to $y \in Y$, then $r$ induces an isomorphism $\pi_z \xrightarrow{\sim} \Gamma_y$ of isotropy subgroups of $\pi$ and $\Gamma$ respectively, as well as an isomorphism of a  neighbourhood of $z$ onto that of $y$, respecting the actions of the isotropy groups. 
 
Since the action of $\pi_o$ is free on $\bh$, by using the invariant direct image functor ${r}^{\pi_o}_{_*}$, the study of $(\Gamma,G)$--bundles on $Y$ reduces to the study of $(\pi,G)$--bundles on ${\bh}$ and thus the study of $(\pi,G)$--bundles on ${\bh}$ reduces to an algebraic problem since $Y$ is a compact Riemann surface and hence a smooth projective curve. 

\brem\label{localdata}By {\em local data} we mean, {\em the subset $\mathcal R \subset X$, together with the ramification indices or signatures $n_i$ and local cyclic coverings of order $n_i$ in formal neighborhoods of the ramification points}.
It is obvious from the above discussion that constructions that involve only the local data are independent of the choice of $Y$, since in principle one could have used the universal cover $\bh$. This could also be seen by using {\em orbifold stacks} constructed from the local data. In the course of this work, we will however work with a fixed $Y$.\erem

\bdefe
A $(\Gamma,G)$--bundle over $Y$ is a principal $G$--bundle $E$ (with a right $G$--action)  together
with a lift of the action of $\Gamma$ on the total space of $E$ as
bundle automorphisms preserving the action of $G$.\edefe

\brem Note that the actions of $G$ and $\Gamma$ on the
total space of $E$ commute which is equivalent to the above condition
that $\Gamma$ acts as automorphisms preserving the action of $G$.\erem

 Now a $(\Gamma,G)$--bundle $E$ on $Y$ is {\em locally} a $(\Gamma_y,G)$--bundle at $y$. Recall that this $(\Gamma_y,G)$--bundle is {\em defined by a  representation}; i.e., if $N_{_y}$ is a sufficiently small $\Gamma_y$--stable formal neighbourhood of $y$, then this bundle is isomorphic to the $(\Gamma_y,G)$--bundle $N_{_y} \times G$, for the {\em twisted}  $\Gamma_{_y}$--action on $E \times G$ given by a representation $\rho_y : \Gamma_{_y} \lr G$, defined as follows: 
\beqa\label{intro2}
\gamma  \cdot (u,g) = (\gamma u, \rho_y (\gamma ) g),~  u \in N_{_y},~ 
\gamma \in \Gamma_{_y}. 
\eeqa 
(See for example Grothendieck \cite[Proposition 1, page 6]{groth}; in the  setting of formal neighbourhoods, see the more recent paper of Teleman-Woodward \cite[Lemma 2.5]{tw}).
\begin{observe} {It is easily seen that these $(\Gamma_y,G)$--bundles given by representations are isomorphic as  $(\Gamma_y,G)$--bundles if and only if the defining representations are equivalent. We call $\rho_y$ the {\em local representations} associated to a $(\Gamma,G)$--bundle}.\end{observe}

\bdefe\label{earlylocaltype} Let $E$ be a $(\Gamma,G)$--bundle on $Y$. The {\em local type} of $E$ at $y$ is defined as the equivalence class of the local representation $\rho_y$ and is denoted by $\tau_y$. \edefe

We denote by $\boldsymbol\tau$, the set $\{\tau_y \mid y \in p^{-1}{\mathcal R}) \}$ (see Definition \ref{localtypeofrep}). Let us denote by 
\beqa\label{intro4}
{Bun_{_Y}^{\boldsymbol \tau}(\Gamma,G)} = \left\{
\begin{array}{l}
\mbox{the set of isomorphism classes of } \\ 
\mbox{$(\Gamma,G)$ bundles with local type $\boldsymbol\tau$}\\ 
\end{array} \right\}
\eeqa

Let $D_{_x} = Spec(A)$. Similarly, for $y \in p^{-1}({\mathcal R})$, let $N_{_y} = Spec(B)$, where $B$ is the integral closure of $A$ in $L = K(\omega)$, where $\omega$ is a primitive $d^{th}$--root of $z$, $d = |\Gamma_y|$ and $z$ is the uniformizer of $A$. Let $p:N_{_y} \to D_{_x} \simeq N_{_y}/{\Gamma_{_y}}$ be the totally ramified covering projection. 
Let $E$ be the $(\Gamma,G)$--bundle on $Y$ and $y \in p^{-1}({\mathcal R})$. Consider the restriction of $E$ to
$N_{_y}$. Then as we have seen above in \eqref{intro2}, as a $(\Gamma_{_y}, G)$ bundle  we can identify $E|_{_{N_{_y}}}$ 
with the trivial bundle $N_{_y} \times G$ together with the twisted $\Gamma_{_y}$--action.

\bdefe\label{unitgrp}  Define ${\sf U}_y$ to be the group:
\beqa
{\sf U}_y = Aut_{_{(\Gamma_{_y},G)}}(E|_{_{N_{_y}}})
\eeqa
of $(\Gamma_{_y}, G)$ automorphisms of $E$ over $N_{_y}$. We call ${\sf U}_y$ {\em the unit group} (or more precisely {\em the local unit group} at $y \in Y$) associated to $E$. \edefe

We work with notations fixed above. Let $\rho_y: \Gamma_{_y} \to G$ be a representation. Let $\ell = rank(G)$ and we represent the maximal torus $T \subset G$  in the diagonal form as follows:
\beqa\label{type1} T = \left[
\begin {array}{cccc}
{t_1}&&&0\\
&.&&\\
&&.&\\

0&&&{t_{\ell}}
\end {array}
\right]\eeqa
where $\{t_1, \ldots, t_{\ell}\}$ is a basis of $X(T)$.
 
Since $\Gamma_{_y}$ is cyclic, we can suppose
that the representation $\rho_y$ of $\Gamma_{_y}$ in $G$ factors through 
$T$ (by a suitable conjugation). 

The action of $\Gamma_{_y}$ on $N_{_y}$ canonically determines a character as follows. The action determines an action of $\Gamma_{_y}$ on the tangent space $T_{_y}$ to $N_{_y}$ at $y$. Since the tangent space to $N_{_y}$ is $1$--dimensional this action is given by a character which we denote by $\chi_{_y}$ (which is of order $d$). Fix a generator $\gamma$ in $\Gamma_{_y}$. The character $\chi_{_y}$ is given by:
\beqa\label{chinought}
\chi_{_y}(\gamma).\omega = \zeta.\omega
\eeqa
where $\zeta$ is a primitive $d^{th}$--root of unity.

\blem\label{initialmumbo} Let $\Gamma_{_y}$ be a cyclic group of order $d$ acting on $N_{_y}$ as above. Then we have a canonical identification
\beqa
Hom(\Gamma_{_y},T) \simeq {\frac{Y(T)}{d. Y(T)}}\eeqa \elem
\begin{proof} This can be seen easily as follows. Observe that $X(\Gamma_y) \simeq {\bz}/{d\bz}$ by the canonical choice of the character $\chi_{_y}$ as in~\eqref{chinought}. Then, we see that
\[
Hom(\Gamma_y,T) = Hom(X(T), X(\Gamma_y)) = Hom(X(T), {\bz}/{d\bz}) = {\frac {Y(T)}{d .Y(T)}}.
\]
\end{proof}
\brem This lemma can be seen in the light of Remark~\ref
{assocbtgrpscheme}. The equivalence class of a representation in $Hom(\Gamma_y,T)$ is given by the conjugacy class of the image of $\gamma$ and hence a point of the Weyl alcove. \erem

We now elaborate this identification  for setting up the notations which play a key role in the next theorem.

Given a representation $\rho_y \in Hom(\Gamma_{_y},T)$, the image $\rho_y(\gamma)$ takes the form
\beqa\label{} \rho_y(\gamma) = \left[
\begin {array}{cccc}
{\chi_{_y}(\gamma)^{a_1}}&&&0\\
&.&&\\
&&.&\\

0&&&{\chi_{_y}(\gamma)^{a_{\ell}}}
\end {array}
\right]\eeqa 
i.e. $\rho_y(\gamma)$ takes the form
\beqa\label{type2} \rho_y(\gamma)= \left[
\begin {array}{cccc}
{\zeta}^{a_1}&&&0\\
&.&&\\
&&.&\\

0&&&{\zeta}^{a_{\ell}}
\end {array}
\right]  {\rm ~with~} a_i \in {\mathbb Z}. \eeqa 

We can suppose that $|a_i| < d$ for all $i$ (or even $0 \leq a_i <d)$ 
and take 
\beqa\label{4'} \eta_i = a_i/d, {\rm ~~so~~ that} ~~ | \eta_i|<1 \eeqa 
Note that the numbers $\{a_1, a_2, \ldots, a_{\ell} \}$ are determined uniquely modulo $d$. Further, this is independent of the choice of $\zeta$.

In terms of the local coordinates $\omega$ and $z$, we may identify  the
function $\omega^{a_i}$ with $z^{\eta_i}$ where $z= \omega^d$.  Define the   
rational  map $\Delta: N_{_y} \lr T$ , or equivalently a morphism
of the punctured disc $N_{_y} - (0)$ as follows:
\beqa\label{4} 
\Delta = \Delta(\omega) = \left[
\begin {array}{cccc}
{\omega}^{a_1}&&&0\\
&.&&\\
&&.&\\

0&&&{\omega}^{a_{\ell}}
\end {array}
\right] 
=
\left[
\begin {array}{cccc}
{z}^{\eta_1}&&&0\\
&.&&\\
&&.&\\

0&&&{z}^{\eta_{\ell}}
\end {array}
\right]\eeqa 

Then we have
\beqa\label{Delta}
\Delta (\gamma u) = \rho (\gamma ) \Delta (u), \quad u \in N_{_y}
\eeqa
where $\Delta$ can be taken as a rational map $\Delta: N_{_y} \lr G$ 
(through $T \hookrightarrow G$). 

Consider the restriction of $\Delta$ to the {\em punctured disc} and view it as a 1PS i.e., $\Delta|_{_{Spec(L)}}:{\bg}_{m,L} \lr G$.  This automatically gives a rational 1--PS of $G$, i.e. an element $ \theta_{\tau_y} \in Y(T) \otimes \bq$ and the key point to note is that 
\beqa\label{thetadelta}
d. \theta_{\tau_y} = \Delta
 ~~i.e.~~ \theta_{\tau_y}  \in {\frac{Y(T)}{d.Y(T)}}\eeqa The association $\rho_y \mapsto \theta_{\tau_y}$ gives explicitly the  identification obtained in Lemma \ref{initialmumbo}. This is precisely what is described in terms of alcoves in Remark~\ref{assocbtgrpscheme}. 
 
\brem\label{choiceofdelta} Note that a choice of $\Delta$ is determined upto (right) multiplication by an element from $G(K)$. \erem

\subsection{\bf The unit group and parahoric groups}   The aim of this section is to prove the following:
\bth\label{weilbt} The unit group ${\sf U}_y$ (Definition \ref{unitgrp})  is isomorphic to  a parahoric subgroup ${\mathcal P}_{_{\theta_{\tau_y} }}(K)$ of $G(K)$ associated to the element $ \theta_{\tau_y}  \in Y(T) \otimes \bq$. Conversely, if ${\mathcal P}_{_\theta}(K)$ is any parahoric subgroup of $G(K)$ then there exists a positive integer $d$, a field extension $L = K(\omega)$ of degree $d$ over $K$ such that
\beqa
{\mathcal P}_{_\theta}(K) \simeq {\sf U}_y
\eeqa
\eeth

\begin{proof} We first give a different description of the elements of ${\sf U}_y$. By \eqref{intro2} a $(\Gamma_{_y}, G)$--bundle on $Y$ gets a $\Gamma_{_y}$--equivariant trivialization; in other words, the $\Gamma_{_y}$--action on $N_{_y} \times G$
is given by a representation $\rho : \Gamma_{_y} \lr G$ 
\beqa\label{eq2.1}
\gamma  \cdot (u,g) = (\gamma u, \rho (\gamma ) g),~  u \in N_{_y},~ 
\gamma \in \Gamma_{_y}. \eeqa
Let  $\phi_0 \in {\sf U}_y$, i.e., the map
\beqa \label{eq2.2}
\phi_0 : N_{_y} \times G \lr N_{_y} \times G. \eeqa
is equivariant for the $\Gamma_{_y}$--action.
Equivariance under $G$ (by right multiplication) implies that 
\[ 
\phi_0(u,g)  =  (u, \phi (u)g)
\]
where $\phi : N_{_y} \lr G$ is a regular map
satisfying the following $\Gamma_y$--equivariance:
\beqa\label{phi} 
\phi (\gamma \cdot u) & = & \rho (\gamma ) \phi (u) \rho (\gamma )^{-1},
u \in N_{_y}, \gamma \in \Gamma_y. \eeqa
We may thus identify ${\sf U}_y$ with the following:
\beqa\label{unitgrp'}
{\sf U}_y = \{\phi:N_{_y} \to G \mid \eqref{phi}~holds \} = Mor^{{\Gamma_y}}(N_y, G)
\eeqa
Since $N_y = Spec(B)$, we can view ${\sf U}_y \subset G(B) \subset G(L)$. 

Let $\Delta$ be as in \eqref{4}. Consider the inner automorphism defined by $\Delta$:
\beqa\label{idelta}
i_{_\Delta}: G(L) \to G(L)
\eeqa
given by $i_{_\Delta}(\eta) = \Delta^{-1}.\eta.\Delta$. Define
\beqa\label{uprime}
{\sf U}'_y := i_{_\Delta}({\sf U}_y)
\eeqa
Let $\psi = i_{_\Delta}(\phi) = \Delta^{-1}.\phi.\Delta$ with $\phi \in {\sf U}_y$. Then we observe that 
\[
\psi (\gamma u) = \psi (u)
\] 
so that $\psi \in G(L)^{\Gamma_{_y}}$. That is, it {\em descends} to a rational function $\tilde{\psi}:D_{_x} \lr G$, where $\tilde{\psi}(z) := \psi(\omega)$. In other words, we get
\beqa
{\sf U}'_y \subset  G(K) = G(L)^{\Gamma_{_y}}.
\eeqa
Note that ${\sf U}'_y$  depends on the choice of $\Delta$ and a different choice of $\Delta$ gives a subgroup which is a conjugate of ${\sf U}'_y$  by an element of $G(K)$ (see Remark \ref{choiceofdelta}).

Then we  {\em claim} the following:
\beqa\label{unitsasparahoric}
{\sf U}'_y = {\mathcal P}_{_{\theta_{\tau_y}}}(K)
\eeqa
where $ \theta_{\tau_y} \in Y(T) \otimes \bq$ is as in \eqref{thetadelta}.
Recall that: 
\beqa
{\mathcal P}_{_{\theta_{\tau_y}}}(K) = \langle T(A),~~~~ U_{r}(z^{m_r(\theta)}A), ~~~r \in R \rangle
\eeqa
Let $\psi \in {\sf U}'_y$ and let $\psi = i_{_\Delta}(\phi)$, with $\phi \in {\sf U}_y$. Thus, 
\[
\phi = \Delta \psi \Delta^{-1}. 
\]
Consider the map $\phi:N_y \to G$. Let $G^o \subset G$ denote the {\em big cell} determined by the roots $R$, (i.e. the inverse image in $G$ of a dense $B$--orbit in $G/B$).

{\em Let us assume for the moment that $\phi(N_y) \in G^o$.}  
In other words, $\phi$ can be described {\em uniquely} as a tuple $\bigl(\{{\phi_r}\}_{_{r \in R}}~ ,~~ \phi_t\bigr)$, with $\phi_r(u) \in U_r$ and $\phi_t(u) \in T$ for $u \in N_y$.

We first consider the tuples $\bigl(\phi_r(u)\bigr)_{r \in R}$ and the corresponding tuple for $\psi$, namely, $ \bigl(\psi_r(u) \bigr)_{r \in R}$, where the ~~$\phi_t~:~N_y \to T$ and $$\{\phi_r, \psi_r: {\bg}_{a,L} \to G \mid r \in R\}.$$ The uniqueness of the decomposition of elements in the big cell and the invariance property of $\phi$ translates into invariance for each of the $\phi_r$ and $\phi_t$. In other words, we have  the following:
\beqa
\phi_r(\omega) = \Delta \psi_r(\omega) \Delta^{-1}. 
\eeqa
i.e.
\beqa 
\phi_r(\omega) = \psi_r(\omega) \omega^{r(\Delta)}
\eeqa
In terms of $\tilde{\psi}$, this gives: 
\beqa\label{psiversusphi}
\phi_r(\omega) = \tilde{\psi_r}(z) z^{\frac{r(\Delta)}{d}}
\eeqa
Now interpreting the condition  that the $\phi$'s are regular functions in the variable $\omega$ at $\omega = 0$, we see that the {\em order of pole} for $\psi_r(z)$ at $z = 0$, is bounded above by $[{\frac{r(\Delta)}{d}}]$ (the biggest integer smaller than or equal to ${\frac{r(\Delta)}{d}}$). In other words $\forall r \in R$,
\beqa
\tilde{\psi_r}(z) \in U_r(z^{-[r(\theta_{_\Delta})]} A) = U_r(z^{m_r(\theta_{_\Delta})} A)
\eeqa
and hence $\tilde{\psi} \in {\mathcal P}_{_{\theta_{\tau_y}}}(K)$.

Now, towards  completing the proof of the claim~\eqref{unitsasparahoric}, since  $\phi_t(u) \in T$, by \eqref{phi} it follows that $\phi_t$ is $\Gamma_y$--invariant and hence, $\tilde{\psi_t} \in T(A)$.

We now take a closer look at the map $\phi:N_y \to G$. In general, the image $\phi(N_y)$ need not be contained in the big cell $G^o$. So we consider the image $\phi(y)$ of the point $y \in N_y$; let $\phi(y) = g_o \in G$. Since the point $y \in N_y$ is $\Gamma_y$--fixed, it implies that $g_o \in G^{\Gamma_y}$. Thus, by \eqref{phi}, the point $g_o$ lies in the centralizer $C_{_G}(\rho(\gamma))$, of $\rho(\gamma)$ in $G$;  the group $C_{_G}(\rho(\gamma))$ is a Levi subgroup $L_{\theta}$ of the standard parabolic subgroup of $G$ determined by the coroot $\theta = \theta_{\tau_y}$. The Levi subgroup can be described in terms of the $u_r:{\bg}_a \to G$ given as in \eqref{root1}, namely $C_{_G}(\rho(\gamma)) = L_{\theta} = \langle T,~~~~ u_{r}(\bc)\mid ~~~r \in R,~~~and~~ m_r(\theta) = (\theta,r) = 0 \rangle$ (see also \ref{stdparahorics}). 

Furthermore, by the equation \eqref{4} which defines the  function $\Delta:{\bg}_m \to T$, it is immediate from \eqref{root1} that ${\Delta}^{-1}.u_r.{\Delta} = u_r$ if $m_r(\theta) = (\theta,r) = 0$. The same obviously holds for the elements of the maximal torus. Hence the elements which commute with $\rho(\gamma)$ also commute with $\Delta$.  This implies immediately that $g_o = i_{_{\Delta}}(g_o)$  and therefore $g_o$ is an element of the parahoric subgroup ${\mathcal P}_{_{\theta_{\tau_y}}}(K)$. 

Now define $\phi_1:N_y \to G$ by $\phi_1(u) = g_o^{-1} \phi(u)$. Then, $\phi_1(y) = 1$ and hence lies in $G^o$. Hence by the openness of $G^o$ and the fact that $N_y$ is a formal neighbourhood of $y$, it follows that $\phi_1(N_y) \subset G^o$. Also, clearly $\phi_1$ satisfies \eqref{phi} and hence by the earlier argument together with the fact that $i_{_{\Delta}}(g_o) \in {\mathcal P}_{_{\theta_{\tau_y}}}(K)$, we see that $i_{_{\Delta}}(\phi) = \psi$ is an element in ${\mathcal P}_{_{\theta_{\tau_y}}}(K)$. This completes the proof of the claim~\eqref{unitsasparahoric} without any assumptions.

Conversely, we show that any parahoric subgroup of $G(K)$ can be identified (upto conjugation by $G(K)$) with a {\em unit group} ${\sf U}_y$. Let $\theta \in {\mathbb E}_{\bq}$ and let ${\mathcal P}_{_{\theta}}(K)$ be a parahoric subgroup. We would like to modify $\theta$ to a $\theta_{\tau_y}$ for a suitable $\Delta \in Y(T)$ (as in \eqref{thetadelta}) so that, interpreted as unit groups we get $ {\mathcal P}_{_{\theta}} (K)  \simeq {\mathcal P}_{_{\theta_{\tau_y}}}(K) \simeq {\sf U}_y$. 

Expressing $\theta$ in terms of generators and clearing denominators, we see that there exists a positive integer $d$ so that $d.\theta \in Y(T)$. Then the obvious choice for $\Delta$ is  simply $d. \theta$, which therefore forces $\Delta \in Y(T)$. The choice of the {\em least} such $d$ makes the choice of the local ramification index canonical.

 Now we view $\Delta$ as a ``rational" map $\Delta:N_{_y} \rightarrow T$ and hence $\Delta$ can be expressed as in \eqref{4}, the $a_i$'s being determined by the following considerations:
for $r \in R$ any root we define
\[
r(\Delta) = d. (\theta,r)
\]

By the discussion following Lemma \ref{initialmumbo}, we have a  $\theta = \theta_{\tau_y} \in {\frac{Y(T)}{d.Y(T)}}$  and the identification of Lemma \ref{initialmumbo} gives the representation $\rho:\Gamma_y \to T \subset G$.  The representation $\rho$ gives the action on the root groups $U_r(B) \subset G(B)$ which are given by (see \ref{root2}):
\beqa\label{morerho}
\rho(\gamma). U_r(B) . \rho(\gamma)^{-1} = U_r(\zeta^{r(\Delta)} B)
\eeqa

Retracing the steps in the first half of the proof, it is easy to see that ${\mathcal P}_{_{\theta}}(K) \simeq {\sf U}_y$, completing the proof of the theorem. \end{proof}

\bnot\label{rhoalpha} Let $\theta \in Y(T) \otimes \bq$. Let $\Delta = d.\theta$ as above. Then we identify $\theta$ with $\theta_{\tau_y}$ and   denote by $\rho_{_\theta}$ the homomorphism $\rho_y: \Gamma_{_y} \to T$ associated to $\theta$ by  Lemma \ref{initialmumbo}.   Note that $\rho_{_\theta}$   acts on the root groups as in~\eqref{morerho}.  
\enot

\brem\label{interiorofalcove}  In the notations used above, if $m_r({\theta_{\tau_y}}) < 1$ for all $r \in R$,   such elements ${\theta_{\tau_y}}$ in ${\mathbb E}_{\bq} = Y(T) \otimes \bq$ are precisely the points in the interior of the alcove $\mathcal A$ (see Remark \ref{interiorofalcove1} and Remark \ref{semisimplealcove}).
\erem

\brem The first half of Theorem \ref{weilbt} can be seen as a consequence of general results on Galois fixed points in Bruhat-Tits buildings and a theorem of Rousseau (\cite[2.6.1]{tits}, \cite{gp}. See also \cite[Section 7]{pr1}). For the converse in Theorem \ref{weilbt} considered in the general setting of Bruhat-Tits theory, we refer the reader to the papers by Gille (\cite[Lemma I.1.3.2]{gille}), Larsen (\cite[Lemma 2.4]{lar}) and Serre (\cite[Proposition 8, p. 546]{ser}). The point of view presented here in terms of unit groups has its origins in the paper of Weil \cite{weil} and Seshadri \cite{pibundles} where completely analogous phenomena are studied in the setting of the general linear group.  The striking fact is that when carried out for semisimple simply connected groups, they yield all parahoric groups when the residue field is of characteristic $0$. \erem

\begin{example}\label{example} Let us now take $G = GL (m)$. We invite the reader to compare this discussion with the one in Weil (\cite[page 56]{weil}). Then we can write 
$\phi = || \phi_{ij} (\omega )||, ~\tilde{\psi} = || \tilde{\psi_{ij}} (z)||,~
1 \leq i,j \leq m$ (as matrices). Then the equation \eqref{psiversusphi} takes the form
\beqa\label{9}\phi_{ij} (\omega ) = \tilde{\psi_{ij}} (z) z^{\alpha_i-\alpha_j}.
\eeqa
We can suppose that $0 \leq \alpha_1 \leq \alpha_2 \leq \cdots \leq \alpha_m < 1$.
Since $|\alpha_i - \alpha_j| < 1$, we deduce easily that 
$\tilde{\psi_{ij}}$ are regular i.e.. ${\sf U}_y \subset
G (A)$.  (To see this suppose that $\tilde{\psi_{ij}}$ is 
{\it not} regular.  Then considered as a function in $\omega$ $(z = \omega^d)$,
$\psi_{ij}$ has a pole of order $\geq d$, whereas $z^{\alpha_i - \alpha_j}$ 
could have only a pole of order $d$ (as a function in $\omega$).  But 
$\phi_{ij}(\omega)$ is regular, which leads to a contradiction).
\end{example}

\brem It is remarked in \cite[Case III, Page 8]{sesramanan}  that it was not clear whether the unit group in the situation considered there is a parahoric subgroup at all. In fact, this is indeed the case as can be seen from Theorem \ref{weilbt}. Moreover, it is not too hard to check by some elementary computations that  the unit group considered in \cite[Case III, Page 8]{sesramanan} does contain the standard Iwahori subgroup but only after a conjugation by a suitable element of $G(K)$ .\erem

\section{The ad\`elic picture of $(\Gamma,G)$--bundles} 
\subsection~ We work with the notations of Section 2. In this section we give a description of $(\Gamma,G)$ bundles analogous to the classical ad\`elic description as in Weil \cite{weil}; it however plays no direct role in the subsequent sections.

Let $E$ be a $(\Gamma,G)$--bundle on $Y$ and let ${Bun_{_Y}^{\boldsymbol \tau }(\Gamma,G)} $ be as in \eqref{intro4}. Since the action of $\Gamma$ on ${Y - p^{-1}({\mathcal R}
)}$ is free, there is a principal $G$--bundle $P$ on $X - {\mathcal R}
$ such that     then $E|_{_{Y - p^{-1}({\mathcal R}
)}} \simeq p^{*}(P)$. Since $G$ is semisimple, by the theorem of Harder \cite{halbein}, $P$ is {\em trivial}. Hence, $E|_{_{Y - p^{-1}({\mathcal R}
)}}$ is also {\em trivial} as a $(\Gamma,G)$--bundle. 

Recall that around each point $y_i \in p^{-1}({\mathcal R}
)$, we have formal neigbourhoods $N_{{y_i}} = Spec~B_{{y_i}}$ with $\Gamma_{{y_i}}$--equivariant trivializations of $E|_{N_{{y_i}}}$ (see \eqref{intro2}). Note that by Beauville-Laszlo (\cite{bl}) any $(\Gamma,G)$--bundle of local type $\boldsymbol\tau$ can be obtained by pathching $E|_{_{Y - p^{-1}({\mathcal R}
)}}$ with the $E|_{N_{{y_i}}}$ 's (see Remark \ref{moregdelta}). 


For simplicity of notation, we assume that ${\mathcal R} = \{x\}$. Two $(\Gamma,G)$ bundles on $Y$ are said to be {\it locally isomorphic at}
$x$ if they are isomorphic as $(\Gamma,G)$--bundles over $p^{-1} (D_{_x}) =
V_1$, $D_{_x}$ a formal neighbourhood of $x$ as above.
We can suppose that $V_1$ is a disjoint union of  $\Gamma_{_y}$--invariant neighbourhood $N_{_y}$ of $y$, $y$ being a point
of $Y$ lying over $x$.  We see that two such bundles are locally isomorphic at $x$ if and only if 
their restrictions to $N_{_y}$ are isomorphic as $\Gamma_{_y}$--bundles. Observe that two  $(\Gamma,G)$--bundles on $Y$ are {\it locally isomorphic at $x$} if they are locally isomorphic in a formal neighborhood of any one point $y \in p^{-1}(x)$. Let 
\beqa
X_1 = X - x, ~~and ~~Y_1 = p^{-1}(X_1)
\eeqa
Recall \eqref{intro2}, the $(\Gamma_y,G)$--bundle $N_{_y} \times G$ given by the twisted  $\Gamma_y$ action given by a representation $\rho_y:\Gamma_y \to G$. 
\beqa\label{eq2.14}
\mbox{\parbox{3.75in}{{\em Let $E_1 := N_{_y} \times G$ with the $(\Gamma_y,G)$--structure given by}
$\gamma  \cdot (u,g) = (\gamma u, \rho_y (\gamma ) g),~  u \in N_{_y},~ 
\gamma \in \Gamma_{_y}.$ }}
\eeqa 
and
\beqa \label{eq2.15}
\mbox{\parbox{3.75in}{{\em Let $E_2 := Y_1 \times G$ with the $(\Gamma,G)$--structure given by}
$\gamma \cdot (u,g) = (\gamma u,g), ~ \gamma \in \Gamma$ and $u \in Y_1$. }}
\eeqa
Thus giving a $(\Gamma,G)$--bundle on $Y$ of local type ${\boldsymbol \tau }$ (see Definition \ref{earlylocaltype}) is giving a
{\em transition function}, i.e  a $(\Gamma_y,G)$--isomorphism:
\beqa \label{eq2.16}
\Theta : E_2 |_{N_y \cap Y_1} \lr E_1 |_{N_y \cap Y_1}.
\eeqa
{\em We denote by $E_{_\Theta}$ the $(\Gamma,G)$--bundle given by the transition function $\Theta$}.

Observe that any transition function $\Theta$ can be viewed as a function $\Delta$ as in \eqref{Delta}. In particular, if we take $\Delta$ as in \eqref{Delta}, then viewed as a {\em transition function} $\Delta$ defines a  $(\Gamma,G)$--bundle, which we denote by $E_{_\Delta}$. We fix this bundle $E_{_\Delta}$ as a {\em base point}.

By Theorem \ref{weilbt}, this choice of $\Delta$ further identifies each unit group ${\sf U}_y'$, $y \in p^{-1}({\mathcal R}
)$  with a parahoric group ${\mathcal P}_{_{{\theta_i}}}(K_{x_i})$, $x_i \in {\mathcal R}$. We fix such an identification.

\bprop\label{settheoretic} For each $x_i \in {\mathcal R}$,  fix a point $y_i \in p^{-1}(x_i)$. Fix further at each $y_i$ local data as in \eqref{eq2.14} and \eqref{eq2.15} and $\boldsymbol\Delta:= \{\Delta_i\}$ as in \eqref{Delta}. Let $K_x$ be the quotient field of the complete local rings $A_x$ at $x \in {\mathcal R}$ and $k[X - {\mathcal R}]$ the ring of functions on the affine curve $X - {\mathcal R}$. Then we have a well-defined set-theoretic identification:
\beqa \label{adelegamma}
 {Bun_{_Y}^{\boldsymbol \tau}(\Gamma,G)}\simeq \Big[\prod_{x \in {\mathcal R}}{\mathcal P}_{_{{\theta_i}}}(K_{x_i}) \backslash ^{\mbox{$\prod_{x \in {\mathcal R}}G(K_x)$}} / G(k[X - {\mathcal R}]) \Big]
\eeqa
where the base point  $E_{_{\boldsymbol\Delta}}$ given by the transition functions $\Delta_i$'s gets identified with the  double coset represented by $1 \in \prod_{x \in {\mathcal R}}G(K_x)$.
\eprop

\begin{proof} For simplicity of notation, we assume again that ${\mathcal R} = \{x\}$. Let $E_{_\Theta}$ and $E_{_\Upsilon}$ be two $(\Gamma,G)$--bundles given by the transition functions $\Theta$ and $\Upsilon$. Then $E_{_\Theta}$ is 
$(\Gamma,G)$--isomorphic to $E_{_\Upsilon}$ if and only if there exist  a $(\Gamma_y,G)$--automorphism $\phi$ of $E_1$ and  
 a  $(\Gamma,G)$--automorphism $\mu$ of $E_2$ such that:
\beqa \label{eq2.18}
\begin{array}{l}
\phi ~\Theta~ \mu = \Upsilon.
\end{array}
\eeqa


We now proceed to give a description of $\phi$ and $\mu$ basing ourselves on the fixed choice of the function $\Delta$.

Observe that by \eqref{eq2.15} the map $\mu$ is given by a morphism:
\beqa \label{eq2.19}
\begin{array}{l}
\mu: Y_1 \times G \lr Y_1 \times G, \\
(u,g) \to (u,\mu(u)g), 
\end{array}
\eeqa
where $\mu(\gamma \cdot u) = \mu (u), \gamma \in \Gamma$. In other words, the map $\mu$ goes down to a morphism $X_1 \lr G$ and we can view $\mu$ as an element in $G(X-x)$.

We now trace the various identifications by restricting the above picture to  $N_{_y}^* = N_{_y} - (0)$; note 
that the $(\Gamma,G)$--isomorphism $\Theta$ is completely characterized by its restriction
to $N_{_y}^*$.  

We observe by \eqref{eq2.15} that the restriction 
of $E_2$ to $N_{_y}^*$ is the $(\Gamma_{_y},G)$--bundle $N_{_y}^* \times G$ over $N_{_y}^*$ with the action
of $\Gamma_{_y}$ given by
\beqa \label{eq2.20}
\begin{array}{c}
\gamma : N_{_y}^* \times G \lr N_{_y}^* \times G, ~~ \gamma \in \Gamma_{_y}\\
\gamma (u,g) = (\gamma u,g).
\end{array}
\eeqa
The restriction of $E_1$ to $N_{_y}^*$ is the $(\Gamma_{_y},G)$--bundle $N_{_y}^* \times G$
on $N_{_y}^*$ with the action of $\Gamma_{_y}$ given by 
\beqa \label{eq2.21}
\begin{array}{l}
\gamma : N_{_y}^* \times G \lr N_{_y}^* \times G \\
\gamma (u,g) = (\rho u, \rho (\gamma )g), ~ \gamma \in \Gamma_{_y}.
\end{array}
\eeqa
The restriction of $\Theta|_{N_{_y}^*}$ of $\Theta$ to $N_{_y}^*$ (denoted again by $\Theta$)  is then a $(\Gamma_{_y},G)$--isomorphism
of the bundle in (\ref{eq2.19}) with the one of (\ref{eq2.18}).
We see easily that $\Theta$ is defined by the map:
\beqa \label{eq2.22}
\begin{array}{c}
N_{_y}^* \times G \lr N_{_y}^* \times G \\
(u,g) \lr (u,\Theta(u)g) \\
\end{array}
\eeqa
where $\Theta : N_{_y}^* \to G$ is such that $\Theta(\gamma \cdot u) = \rho (\gamma ) \Theta (u)$.

Recall that the map $\Delta$ as in (\ref{4}) is a morphism $N_{_y}^* \lr G$
and has similar properties.  Thus we can write
\beqa \label{eq2.23}
\begin{array}{l}
\Theta = \Delta \Theta{_o} ~~{\rm such ~that}~~~ 
\Theta{_o}(\gamma u) = \Theta{_o} (u) 
\end{array}
\eeqa
i.e.. $\Theta{_o}$ descends to a regular map $D_{_x}^* \lr G$, $D_{_x}^* = D_{_x} - (0)$.

The equivalence relation (\ref{eq2.18}) therefore takes the following form:
\beqa \label{eq2.25}
\phi ~(\Delta \Theta{_o})~ \mu =  \Upsilon
\eeqa
Multiplying on either side by $\Delta^{-1}$ we get
\beqa\label{eq2.251} 
({\Delta^{-1} \phi \Delta})~ \Theta{_o}~ \mu =  \Delta^{-1} \Upsilon = \Upsilon_{o}.
\eeqa
By the proof of Theorem \ref{weilbt},  $\phi$
identifies with an element $\psi (= i_{\Delta}(\phi))$ of the unit group ${\sf U}_y'$ and we can write \eqref{eq2.251} as
\beqa\label{eq2.252} 
\psi~ \Theta{_o}~ \mu =  \Delta^{-1} \Upsilon = \Upsilon_{o}.
\eeqa

 Therefore,  $\Theta{_o}
\in G(K_x)$ and $\psi \in {\sf U}_y' = {\mathcal P}_{_{{\theta}}}(K_{x})$  and  by \eqref{eq2.19}, $\mu$
becomes a regular map $X_1 \lr G$ i.e.. $\mu \in G(X - x)$. 

From (\ref{eq2.252}), we conclude that $\Theta$ and $\Upsilon$ give isomorphic $(\Gamma,G)$--bundles if and only if $\Theta_o$  and $\Upsilon_o$ are equivalent by the double coset relation, i.e give the same point in the double coset space which we denote by $[\Theta_o]$. If $\Theta_o \in G(K_x)$ gives $[\Theta_o]$ in  the double coset space, by using \eqref{eq2.23} and the choice of $\Delta$, we can reverse the process to get $\Theta$ and hence $E_{_\Theta}$. Thus, we get the following set-theoretic identification:
\beqa \label{eq2.26}
{Bun_{_Y}^{\boldsymbol \tau }(\Gamma,G)}  \simeq \Big[{\mathcal P}_{_{{\theta}}}(K_{x}) \backslash ^{\mbox{$G(K_x)$}} / G(X - x) \Big]
\eeqa
where $E_{_\Theta} \mapsto [\Theta_o]$. 

Note that the base point $E_{_{\Delta}}$ given by the transition functions $\Delta$ gets mapped to the {\em identity coset}, i.e. represented by $1 \in G(K_x)$ in the double coset space.\end{proof}

\section{Invariant direct image}
\subsection ~In this section we study the torsor-analogue of the sheaf theoretic invariant direct image defined by Grothendieck \cite{groth}. The remarks in this section owe much to key inputs from Brian Conrad and Pramathanath Sastry.

\bdefe\label{galoiscover} Let $p:W \to T$ be a finite flat surjective morphism of normal, integral noetherian schemes such that the field extension of the function fields is Galois with Galois group $\Gamma:= {\text {Gal}}(k(W)/k(T))$. Observe that $\Gamma$ acts on $W$ as $T$--automorphisms and $T = W/{\Gamma}$. Such a morphism $p:W \to T$ is called a {\em Galois covering} with Galois group $\Gamma$.\edefe


Following \cite{blr}, we can define the {\em direct image functor} $p_*$ as the Weil restriction of scalars, i.e., we have a group functor $p_{_*}(\eG) := {\text{Res}}_{_{W/T}}(\eG)$ with the following property; for any $T$--scheme $S$, we have a canonical bijection, i.e, the {\em adjunction}:
\beqa
Hom_T(S, p_{_*}(\eG) \simeq Hom_W(S \times_T W, \eG)
\eeqa
which is functorial in $S$ and $\eG$. 

We assume that $\eG$ is an {\em affine group scheme} over $W$, so that $p_{_*}(\eG)$ is  representable by a  group scheme (see \cite[Theorem 4 and Proposition 6]{blr}). Suppose also that the $\Gamma$--action lifts to an action on the group scheme $\eG$, in such a manner that the `multiplication map' and the `inverse map' on $\eG$ are equivariant. We will term such a group scheme a {\em $\Gamma$--group scheme} on $W$.  

Let $S$ be a scheme over $T$, and $f \in p_{_*}(\eG)(S) = Hom_{_{W}}(S \times_T W,\eG)$, and let $\gamma \in \Gamma$. 

There is a left action of $\Gamma$ on $\eG$ and  a left action on $S \times_T W$ induced by its action on $W$. 
 This induces a natural right action of $\Gamma$ on $p_{_*}(\eG)(S)$ given by:
\beqa
(f.\gamma)([s,w]) := \gamma^{-1}.f({\gamma}.[s,w]), ~~~~[s,w] \in S \times_T W
\eeqa
We can now take the {\em fixed point} subscheme under the action of $\Gamma$. The general results on fixed point subschemes given in \cite[Section 3]{bass} can be applied to our situation {\em since we are in characteristic $0$} and we get a canonically defined {\em smooth} closed $X$--subgroup scheme $p_{_*}(\eG)^{\Gamma} \subset p_{_*}(\eG)$. This is representable in our case since $p_{_*}(\eG)$ is representable.
\bdefe\label{invariantdirectimage}(Invariant direct image) Let $p:W \to T$ be as above and let $\Gamma = Gal(W/T)$. Let $\eG$ be a smooth affine $\Gamma$--group scheme over $W$. We define the {\em invariant direct image} of $\eG$ to be:
\beqa
p^{\Gamma}_{_*}(\eG):= p_{_*}(\eG)^{\Gamma}
\eeqa
i.e., for any $T$--scheme $S$, we have $p^{\Gamma}_{_*}(\eG)(S) = \eG(S \times_T W)^{\Gamma}$. 

More generally, let $E$ be any affine scheme over $W$ with a lift of the $\Gamma$ action. Then we define the {\em invariant direct image}  of $E$ to be $p^{\Gamma}_{_*}(E):= p_{_*}(E)^{\Gamma}$.
\edefe

\blem\label{a} Let $p:W \to T$ be a finite flat surjective morphism of noetherien schemes. Let  $\eG$ be a smooth affine group scheme on $W$ and $E$ a $\eG$--torsor on $W$.  Then $p_{_*}(E)$ is a $p_{_*}(\eG)$--torsor on $T$.\elem

\begin{proof} The hypothesis implies that $p_{_*}(\eG)$ and $p_{_*}(E)$ exist as smooth schemes over $T$. The lemma follows immediately from the property that the direct image functor $p_{_*}$  respects fibre products (this is immediate from the functorial definition of restriction of scalars, see for example \cite[Proposition A.5.2]{conrad}). Applying $p_{_*}$ to the action map $\eG \times_{_W} E \to E$ it gives $p_{_*}(\eG) \times_{_T} p_{_*}(E) \to p_{_*}(E)$. Moreover, we have an isomorphism:
\beqa
\eG \times_{_W} E \simeq E \times_{_W}E
\eeqa
Now again apply $p_{_*}$ to get the desired isomorphism 
\beqa\label{fromcgp}
p_{_*}(\eG) \times_{_T} p_{_*}(E) \simeq p_{_*}(E) \times_{_T} p_{_*}(E)
\eeqa
This is also given in \cite[Corollary A.5.4(3)]{conrad}  but with a more complicated proof.
\end{proof}
If the $\Gamma$--action lifts to an action on a $W$--group scheme so that the `multiplication map' and the `inverse map'  are equivariant, then we will term such a group scheme a {\em $\Gamma$--group scheme} on $W$.  
\blem\label{b} Suppose further that $p:W \to T$ is a Galois cover with Galois group $\Gamma$ (Definition \ref{galoiscover}). Let $\eG$ be a $\Gamma$--group scheme on $W$ and $E$ a $(\Gamma,\eG)$--torsor on $W$. Let $p^{\Gamma}_{_*}(\eG) = \eH$ and $p^{\Gamma}_{_*}(E) = F$. Then, $F$ is a $\eH$--torsor.\elem 
\begin{proof} For the first part, apply  the {\em fixed point} functor to \eqref{fromcgp} i.e., 
\beqa
(p_{_*}(\eG) \times_{_T} p_{_*}(E))^{\Gamma}  \simeq (p_{_*}(E) \times_{_T} p_{_*}(E))^{\Gamma} 
\eeqa
which gives
\beqa
\eH \times_{_T} F \simeq F \times_{_T} F
\eeqa
proving that $F$ is a $\eH$--torsor on $X$. The smoothness of $F$ over $T$ holds as well since we work in char $0$ (cf. \cite[Section 3]{bass}).
\end{proof}

\bth\label{invariantdi} Let ${\sf Bun}_{_W}(\Gamma,\eG)$ and ${\sf Bun}_{_T}(\eH)$ denote the  stacks of $(\Gamma,\eG)$--torsors on $W$ and $\eH$--torsors on $T$ respectively. Then the functor
\beqa
p^{\Gamma}_{_*}:{\sf Bun}_{_W}(\Gamma,\eG) \to {\sf Bun}_{_T}(\eH)
\eeqa
is an isomorphism of stacks.
\eeth
\begin{proof} Lemma \ref{b}  shows that $p^{\Gamma}_{_*}$ gives a functor between the stacks. We now construct the candidate for the inverse.

Observe that the inclusion $p_{_*}^{\Gamma}(\eG) = \eH \hra p_{_*}(\eG)$ gives by ``adjunction" the morphism:
\beqa\label{c}
p^*(\eH) \to \eG
\eeqa
of group schemes over $W$. Let $F \in {\sf Bun}_{_T}(\eH)(S)$ be a $\eH$--torsor on a $T$--scheme $S$. Let $p \times Id_S = q:W \times_T S \to S$ be the induced morphism. Observe that $q^*(F)$ becomes a ${q^*(\eH)}$--torsor on $W$ and via \eqref{c} we get the associated $\eG$--torsor
\beqa
q^*(F) \times ^{q^*(\eH)} \eG
\eeqa
We observe that the $\eG$--torsor $q^*(F) \times ^{q^*(\eH)} \eG$ is a $(\Gamma,\eG)$--torsor, where the $\Gamma$--action comes from the underlying $\Gamma$--action on $\eG$. 
We also get the natural $\Gamma$--equivariant morphism:
\beqa
q^*(F) \to q^*(F) \times ^{q^*(\eH)} \eG
\eeqa
Now by pushing down this morphism using $q^{\Gamma}_{_*}$, we get
\beqa
F \to q^{\Gamma}_{_*}(q^*(F) \times ^{q^*(\eH)} \eG)
\eeqa
To check this last map is an isomorphism, we can restrict to \'etale neighbourhoods on $T$ where $F$ is trivial (i.e, isomorphic to $\eH$ as an $\eH$--torsor); but this is obvious. This shows that $q^{\Gamma}_{_*}(q^*(F) \times ^{q^*(\eH)} \eG) \simeq F$.

We need to check  that this construction provides an equivalence of categories.

Suppose that $E$ is a $(\Gamma,\eG)$--torsor on $W$ and $F = p^{\Gamma}_{_*}(E)$ is a $\eH$--torsor. Now $p^*(F)$ is a $p^*(\eH)$--torsor. Therefore via \eqref{c}, taking associated constructions we get 
a $\eG$--torsor $p^*(F) \times ^{p^*(\eH)} \eG$.

Again by adjunction applied to the inclusion $F = p^{\Gamma}_{_*}(E) \hra 
p_{_*}(E)$, we get the morphism
\beqa
p^*(F) \to E
\eeqa
and hence a morphism
\beqa\label{d}
p^*(F) \times ^{p^*(\eH)} \eG \to E
\eeqa
of $\eG$--torsors.

{\it Claim}: The morphism \eqref{d} is an isomorphism of $(\Gamma,\eG)$--torsors. 

{\it Proof of Claim}:  Since the map $p$ is finite, a cofinal system of \'etale neighbourhoods of the fiber $p^{-1}(t)$ of a point $t \in T$ is given by pullbacks of \'etale neighbourhoods of $t \in T$. This is a consequence of the compatibility of formation of strict henselization with respect to finite base change of algebras (cf. EGA IV.4, 18.8.10). This allows us to work \'etale locally on $T$. 

Thus we may assume that $F$ is trivial on $T$ and by the discussion above we may assume that $E$ is also a trivial $\eG$--torsor on $W$. This reduces to verifying that the map $p^*(\eH) \times ^{p^*(\eH)} \eG \to \eG$ is an isomorphism of trivial $\eG$--torsors, which is obvious. This completes the proof of the theorem.\end{proof}

\brem As the referee pointed out to us, the equivalence can also be deduced from  the observation that both the stacks in Theorem \ref{invariantdi} are gerbes. \erem

\brem The notion of invariant direct image using Weil restriction of scalars is implicit in Edixhoven \cite{bass} and also in Pappas-Rapoport \cite{pr1}.\erem

\brem  Let ${\mathcal O}_W(\eG)$ be the sheaf of groups on $W$ for the \'etale topology associated to the group scheme $\eG$. In fact, it is a $\Gamma$--sheaf of groups. The fact used in the argument above, namely \'etale trivializing neighbourhoods in $W$ can be chosen as inverse images of \'etale opens from $T$ shows firstly that any $\eG$--torsor can be trivialized in such \'etale neighbourhoods. In particular, taking $(\Gamma, \eG)$--torsors, this immediately gives a natural isomorphism of cohomology sets: 
\beqa
H_{_{\acute{e}t}}^1(W,\Gamma,{\mathcal O}_W(\eG)) \simeq H_{_{\acute{e}t}}^1(T,{\mathcal O}_T(p^{\Gamma}_{_*}(\eG)) 
\eeqa
This identifies the isomorphism classes of $(\Gamma,\eG)$--torsors on $W$ with isomorphism classes of $p^{\Gamma}_{_*}(\eG)$--torsors on $T$. The proof given above for the theorem gives a canonical identification and hence a stronger statement on stacks.
\erem
\section{Bruhat-Tits group schemes and torsors}

\subsection{A $\Gamma$--group scheme on $Y$} We now revert to the notations in Section 2, where $p:Y \to X$ ramified over $\mathcal R$. The following construction plays an important role in the subsequent sections. 

\begin{notation}\label{bigdelta}  Fix a $(\Gamma,G)$--bundle $F$ of local type $\boldsymbol\tau$; let 
\beqa
{\mathcal G}_{_F} := F \times^{G} G
\eeqa
denote the associated ``adjoint" group scheme associated to $F$, $G$ acting on itself by inner conjugation.
\end{notation}

Recall that any $(\Gamma,G)$--bundle of local type $\boldsymbol\tau$ is locally isomorphic to any preassigned  $(\Gamma,G)$--bundle of local type, in particular to the fixed bundle $F$. This therefore gives an identification of  $(\Gamma,G)$--bundles of type $\boldsymbol\tau$ with those that are locally modelled after the fixed bundle $F$. We give a formal shape  to this intuitive picture (see Grothendieck \cite[Proposition 4.5.2]{kansas}).

The bundle $F$ can be viewed as a left ${\mathcal G}_{_F}$--torsor, where the action is by automorphisms. In the sense of Giraud \cite{giraud},  $F$  is  a $({\mathcal G}_{_F},G)$--bitorsor.

For any $G$--torsor $E$   on $Y$  coming with a right $G$--action, let $E^{op}$ be the $G$--torsor with the induced left action: $g.x := x g^{-1}$. We then have the ``contracted product":
\beqa
E \wedge^G F^{op} := \frac{E \times_Y F}{(xg,y) \sim (x,g.y)}
\eeqa
It is a fact that the sheaf of local sections of $E \wedge^G F^{op}$ is the sheaf $Isom(E,F)$ of local isomorphisms of $E$ with $F$. Since we work with affine group schemes, by usual descent for affine schemes, the contracted product is representable as a scheme. 

Since $F$ is a $({\mathcal G}_{_F},G)$--bitorsor, it follows that $F^{op}$ is a $(G,{\mathcal G}_{_F})$--bitorsor. Thus the contracted product  $E \wedge^G F^{op}$ is  a right ${\mathcal G}_{_F}$--torsor on $Y$; in fact we get the equivalence of stacks:
\beqa
{\sf Bun}_{_Y}(G) \simeq {\sf Bun}_{_Y}({\mathcal G}_{_F})
\eeqa
given by $E \mapsto E \wedge^G F^{op}$

Now let $E$ be a $(\Gamma,G)$--bundle of local type $\boldsymbol\tau$.
Hence $E$ is locally $\Gamma$--isomorphic to $F$. Since ${\mathcal G}_{_F}$ is a $\Gamma$--group scheme, the association  $E \mapsto E \wedge^G F^{op}$ in fact induces an identification:
\beqa\label{firststep}
{{\sf Bun}_{_Y}^{\boldsymbol \tau }(\Gamma,G)} \simeq {{\sf Bun}_{_Y}(\Gamma,{\mathcal G}_{_F})}
\eeqa
the identification being obviously dependent on the choice of $F$.

\bsem\label{localbtgrpscheme}We now return to the setting in Section 2, i.e. $p:N_y \to D_x$. Recall that $\Gamma_y = Gal(N_y/D_x)$. 
Let $F_y$ be any $(\Gamma_y,G)$--bundle of local type $\tau_y$ and therefore given as in \eqref{intro2}. Since the underlying $G$--bundle is trivial,  the associated adjoint group scheme ${\mathcal G}_{_{F_y}}$ is isomorphic to the product $G \times N_y$. Hence the sections over $N_y = Spec(B)$ are given by ${\mathcal G}_{_{F_y}}(B) \simeq G(B)$. As has been observed in \eqref{unitgrp'}, the local unit group of $(\Gamma_y,G)$--automorphisms ${\sf U}_y$ is a subgroup of $G(B)$. \esem

The content of the first half of Theorem~\ref{weilbt} is that
\beqa\label{parahoricviainvariants}
{\mathcal G}_{_{F_y}}(B)^{\Gamma_y} \simeq {\mathcal P}_{_{\theta_{\tau_y}}}(K)
\eeqa

\bprop\label{morphicequality} Let  ${\mathcal G}_{_{\theta_{\tau_y}}}$ be a Bruhat-Tits group scheme defined by the parahoric group ${\mathcal P}_{_{\theta_{\tau_y}}}(K)$. Let $D_x = Spec~A$ and $N_y = Spec~B$. Let ${\mathcal G}_{_{{F_y}}}$ be the $\Gamma$--group scheme on the fixed $(\Gamma_y,G)$--bundle ${F_y}$. Then
\beqa
{\mathcal G}_{_{\theta_{\tau_y}}} \simeq p^{\Gamma_y}_{_*}({\mathcal G}_{_{{F_y}}})
\eeqa
In particular, if ${\mathcal G}_{_{\theta}}$ is any Bruhat-Tits group scheme on $Spec~A$, then by choosing $\theta_{\tau_y}$ suitably, we can realize ${\mathcal G}_{_{\theta}}$ as $p^{\Gamma_y}_{_*}({\mathcal G}_{_{{F_y}}})
$; for a scheme $S$ over $\bc$, we have the identification:
\beqa\label{morphicequality1}
Mor^{\Gamma_y}(N_y \times S, G) \simeq Mor(D_x \times S,{\mathcal G}_{_{\theta}})
\eeqa
\eprop

\begin{proof} Firstly, since ${\mathcal G}_{_{{F_y}}}$ is an affine group scheme it follows that $p_{_*}({\mathcal G}_{_{{F_y}}})$ and $p^{\Gamma_y}_{_*}({\mathcal G}_{_{{F_y}}})
$ are both representable as affine group schemes over $A$ (see \cite[Theorem 4 and Proposition 6]{blr}).

By Bruhat-Tits (\cite[Section 1.7]{bt}), the smooth group scheme ${\mathcal G}_{_{\theta_{\tau_y}}}$  on $D_x$ is uniquely determined by its $A$--valued points which is the parahoric group ${\mathcal P}_{_{\theta_{\tau_y}}}(K)$.  

By the functorial property of the functor $p^{\Gamma_y}_{_*}$, we see (by \eqref{}) that 
\beqa
p^{\Gamma_y}_{_*}({\mathcal G}_{_{{F_y}}})(A) = p_{_*}({\mathcal G}_{_{{F_y}}})^{\Gamma_y}
= {\mathcal G}_{_{{F_y}}}(B)^{\Gamma_y} \simeq {\mathcal P}_{_{\theta_{\tau_y}}}(K)
\eeqa
Thus, by the uniqueness of the Bruhat-Tits group scheme,  we have an isomorphism of $Spec(A)$--group schemes: ${\mathcal G}_{_{\theta_{\tau_y}}} \simeq p^{\Gamma_y}_{_*}({\mathcal G}_{_{{F_y}}})$.

The identification \eqref{morphicequality1} now follows from the functorial properties of restriction of scalars and fixed point schemes since, 
\beqa
Mor^{\Gamma_y}(N_y \times S, G) = p^{\Gamma_y}_{_*}({\mathcal G}_{_{{F_y}}})(D_x \times S) = {\mathcal G}_{_{\theta}}(D_x \times S)
\eeqa 

The fact that any Bruhat-Tits group scheme can be realized this way follows from the converse in Theorem \ref{weilbt}.
\end{proof}

\subsection{\bf Bruhat-Tits group schemes and patching} 
By the main theorem of Bruhat-Tits (\cite{bt}), there exist smooth group schemes ${\mathcal G}_{_\Omega}$ over $Spec(A)$ such that the group ${\mathcal G}_{_\Omega}(A) =  {\mathcal P}_{_\Omega}(K)$. 

\bdefe\label{globalbt} A smooth group scheme $\mathcal G$ over $X$ is called a parahoric Bruhat-Tits group scheme if there is a finite subset ${\mathcal R} = \{x_i\}$ of $X$, and formal neigbourhoods ${D_{x_i}}$ at the $x_i$  together with a collection of subset $\Theta \subset {\mathbb E}^m$ such that 
\beqa{\mathcal G}|_{_{_{X - {\mathcal R}}}} \simeq G \times (X - {\mathcal R}), ~~~~~~ {\mathcal G}|_{D_{x_i}} \simeq {\mathcal G}_{_{\Theta_i}}, x_i \in {\mathcal R}.\eeqa  
\edefe

If $\Omega = \{\Omega_i\}$ is a collection of facets then we denote a parahoric Bruhat-Tits group scheme defined by local group schemes ${\mathcal G}_{_{\Omega_i}}$ by ${\mathcal G}_{_{\Omega,X}}$.  If $\boldsymbol\theta = \{\theta_i\} \in (Y(T) \otimes \bq)^m$ are chosen in the interior of the facets  $\Omega_i$, then  we have an isomorphism ${\mathcal G}_{_{\Omega,X}} \simeq {\mathcal G}_{_{\boldsymbol\theta,X}}$

Conversely, given local Bruhat-Tits groups schemes ${\mathcal G}_{_{\Theta_i}}$, one can construct a parahoric Bruhat-Tits group scheme using the following  patching result from \cite[Lemma 3.18]{cgp} attributed to Raghunathan and Ramanathan:

\blem\label{gettingbt} Let $X$ be a smooth projective curve and $k(X)$ be its function field. Let $x \in X$ and let $A_{x}$ be the completion of ${\mathcal O}_{X,{x}}$ and $K_x$ the completion of $k(X)$. Assume that we are given a triple $(G_1, G_2, f)$ consisting of:
\begin{enumerate}
\renewcommand{\labelenumi}{(\alph{enumi})}
\item An affine group scheme $G_1$ over $U = X - x$ of finite type.
\item An affine and finitely presented group scheme $G_2$ over $A_x$.
\item A $K_x$--group scheme isomorphism $f:G_1 \times_U {K_x} \simeq G_2 \times_U {K_x}$.
\end{enumerate}
Then there exists a group scheme $\mathcal G$, affine and of finite type over $X$ such that ${\mathcal G} \times_X U \simeq G_1$ and ${\mathcal G} \times_X {A_x} \simeq G_2$ and both isomorphisms are compatible with $f$. Furthermore, if $G_i$ are smooth then so is ${\mathcal G}$. \elem

\brem\label{moregdelta1} The gluing result of  Beauville-Laszlo (\cite{bl}) more general than the lemma above, shows  that  any $\mathcal G$--torsor $E$ on $X$ can be obtained by gluing the trivial torsor on some open subset $U \subset X$ and the trivial torsors on the formal completions at the points ${\mathcal R} = X - U$. Similarly,  any $(\Gamma,G)$--bundle of local type $\boldsymbol\tau$ on $Y$ is obtained by patching  
as was explained in the beginning of Section 3. \erem 

\brem The parahoric Bruhat-Tits group schemes defined above are a little more restrictive than the ones defined by Pappas and Rapoport \cite{pr1}; they do not make the assumption that the group schemes are generically split. \erem

\brem Observe firstly that in Lemma \ref{gettingbt} we can take a finite set of points $x_i \in X$ for the patching; it follows that given a finite ${\mathcal R}
 \subset X$,  and a collection of subset $\Theta \subset {\mathbb E}^m$ together with {\em patching data} $f = \{f_i\}_{i=1}^{m}$ as in Lemma above,  we have a parahoric Bruhat-Tits group scheme ${\mathcal G}_{_{\Theta,X}}$ with  ${\mathcal R}$ being the points of ramification of ${\mathcal G}_{_{\Theta,X}}$.\erem

\brem\label{moregdelta} Let $F$ be a fixed $(\Gamma,G)$--bundle of local type $\boldsymbol\tau$. The group scheme ${\mathcal G}_{_{F}}$ constructed in Notation~\ref{bigdelta} can be viewed as
one obtained by gluing the local group schemes ${\mathcal G}_{_{F_{y}}}$'s on $\{N_y\}_{_{y \in p^{-1}({\mathcal R})}}$ (see \ref{localbtgrpscheme}) along with the constant group scheme $G \times (Y - p^{-1}({\mathcal R}))$, the patching data coming from the transition functions of the bundle $F$.\erem

\bth\label{morphicmore}~Let $F$ be a fixed $(\Gamma,G)$--bundle of local type $\boldsymbol\tau$ on $Y$. Let ${\boldsymbol\theta_{\boldsymbol\tau}} = \{\theta_i\} \in (Y(T) \otimes \bq)^m$ be the point associated to $\boldsymbol\tau$. Then  the invariant direct image $p^{\Gamma}_{_*}({\mathcal G}_{_{F}})$ is a parahoric Bruhat-Tits group scheme of the form ${\mathcal G}_{_{\boldsymbol\theta_{\boldsymbol\tau},X}}$. Conversely, let ${\mathcal G}_{_{\boldsymbol\theta_{\boldsymbol\tau},X}}$ be any parahoric Bruhat-Tits group scheme on $X$, with $\mathcal R$ its set of ramifications. Then, there exists a Galois cover $p:Y \to X$ with Galois group $\Gamma$, and a $(\Gamma,G)$--bundle $F$ of local type $\boldsymbol\tau$  with its adjoint $\Gamma$--group scheme ${\mathcal G}_{_{F}}$ on $Y$,  such that $p^{\Gamma}_{_*}({\mathcal G}_{_{F}}) \simeq {\mathcal G}_{_{\boldsymbol\theta_{\boldsymbol\tau},X}}$.

\eeth

\begin{proof} Since the group scheme 
${\mathcal G}_{_{F}}$ is affine over $Y$, $p^{\Gamma_y}_{_*}({\mathcal G}_{_{F}})$ is representable as a smooth affine group scheme over $X$ (see \cite[Theorem 4 and Proposition 6]{blr}). 

Since the action of $\Gamma$ on ${Y - p^{-1}({\mathcal R}
)}$ is free, there is a principal $G$--bundle $P$ on $X - {\mathcal R}
$ such that     then $F|_{_{Y - p^{-1}({\mathcal R}
)}} \simeq p^{*}(P)$. Since $G$ is semisimple, by the theorem of Harder \cite{halbein}, $P$ is {\em trivial}. 
Therefore, the $(\Gamma,G)$--bundle $F$ when restricted to $Y- p^{-1}(\mathcal R)$ is trivial as a $(\Gamma,G)$--bundle. Hence ${\mathcal G}_{_{F}}$ is the split group scheme over $Y - p^{-1}(\mathcal R)$. The result now follows from Proposition \ref{morphicequality}, the patching data $f$ being the one pushed down from that of ${\mathcal G}_{_{F}}$.

For the converse, observe that  locally, the statement in the corollary is simply the converse in Proposition  \ref{morphicequality}. The global statement now follows since the patching data $f$ gives the gluing needed in Remark \ref{moregdelta1} which gives the recipe to construct $F$ globally. 
\end{proof}

\brem An interesting consequence of Theorem \ref{morphicmore} is that any parahoric Bruhat-Tits group scheme  which is {\em generically split} is isomorphic to the invariant direct image of a group scheme ${\mathcal G}_{_{F}}$ for a choice of $(\Gamma,G)$--bundle $F$. Moreover, this characterizes such group schemes. Observe that in the patching Lemma \ref{gettingbt} one need not assume that the group scheme is generically split. Using this one can show that the parahoric Bruhat-Tits group schemes considered by Pappas and Rapoport  can also be realized as invariant direct images of $\Gamma$--group schemes, which however need not be of the form ${\mathcal G}_{_{F}}$ for a $(\Gamma,G)$--bundle $F$ of type $\boldsymbol\tau$.\erem
  
\subsection{Torsors under Bruhat-Tits group schemes} 

Let ${\mathcal G}_{_{{\boldsymbol\theta},X}}$ be a Bruhat-Tits group scheme given by the local data ${\boldsymbol\theta} \in (Y(T) \otimes \bq)^m$.
 
Let ${\sf Bun}_{_Y}^{\boldsymbol \tau }(\Gamma,G)$ and ${\sf Bun}_{_X}({\mathcal G}_{_{{\boldsymbol\theta},X}})$ be the moduli stacks of $(\Gamma,G)$--bundles of type ${\boldsymbol \tau }$  on $Y$ and of ${\mathcal G}_{_{{\boldsymbol\theta},X}}$--torsors on $X$ respectively. We now have the following key theorem: 

\bth\label{stackygammaversusparahoric} Let  ${\mathcal G}_{_{\boldsymbol\theta_{\boldsymbol\tau},X}}$ be  as above. Let $p:Y \to X$  be as in Theorem \ref{morphicmore}. Then  the stack  ${\sf Bun}_{_Y}^{\boldsymbol \tau }(\Gamma,G)$ is isomorphic to the stack ${\sf Bun}_{_X}({\mathcal G}_{_{\boldsymbol\theta_{\boldsymbol\tau},X}})$.
\eeth
\begin{proof} By Theorem \ref{morphicmore}, there exist a $(\Gamma,G)$--bundle $F$ of local type $\boldsymbol\tau$ on $Y$ such that $p^{\Gamma_y}_{_*}({\mathcal G}_{_{F}}) \simeq {\mathcal G}_{_{\boldsymbol\theta_{\boldsymbol\tau},X}}$ . By \eqref{firststep} we have an isomorphism ${\sf Bun}_{_Y}^{\boldsymbol \tau }(\Gamma,G) \simeq {\sf Bun}_{_Y}(\Gamma,{\mathcal G}_{_{F}})$. By Theorem \ref{morphicmore} and Theorem \ref{invariantdi}, we get the isomorphism ${\sf Bun}_{_Y}(\Gamma,{\mathcal G}_{_{F}})  \xrightarrow{p^{\Gamma}_{_*}} {\sf Bun}_{_X}({\mathcal G}_{_{\boldsymbol\theta_{\boldsymbol\tau},X}})$. This proves the theorem. \end{proof}

 
\brem\label{smallgenus} From the arguments and the results in the preceding pages, it would be clear to the reader that Theorem \ref{morphicmore} and Theorem \ref{stackygammaversusparahoric} also hold when $X$ is $\bp^1$ or an elliptic curve. Note however that $\Gamma$-covers $Y$ will exist only when we assume that $|\mathcal R| \geq 3$ for $X = \bp^1$, or $\mathcal R \neq \emptyset$, when $X$ is an elliptic curve. This is so since the upper half space $\bh$ is the universal ramified cover with the given signature even in these cases.\erem

\brem This theorem is the exact analogue of the fact that the {\em invariant direct image functor} $p^{\Gamma}_{_\ast}$  sets up an isomorphism between the functor of $\Gamma$--vector bundles and that of parabolic vector bundles; this is precisely the point of view in Grothendieck\cite{groth}, Seshadri \cite{pibundles} and Mehta-Seshadri~\cite{ms}.\erem

\section{Stability and semistability} 
The aim of this section is to introduce the notion of semistability and stability of torsors under parahoric Bruhat-Tits group schemes introduced in the last section.

\subsection{\bf Parahoric torsors}   Let ${\mathcal G}_{_{\Omega,X}}$ be as before a Bruhat-Tits group scheme on the curve $X$ associated to a collection of {\em facets}  $\Omega = \{\Omega_i\}_{i=1}^{m}$, with $|\mathcal R| = m$. 

\bdefe\label{quasiparahoric}  A {\bf quasi-parahoric} torsor $\eE$ is a ${\mathcal G}_{_{\Omega,X}}$--torsor on $X$. \edefe


\bdefe\label{parahorictorsor}  A {\bf parahoric torsor} is  a pair $(\eE, {\boldsymbol\theta})$ consisting of
\begin{enumerate}
\item A ${\mathcal G}_{_{\Omega,X}}$--torsor $\eE$, i.e. a quasi-parahoric torsor on $X$ and  
\item Weights, i.e elements ${\boldsymbol\theta} = \{\theta_i\} \in (Y(T) \otimes \bq)^m$ in the interior of the facets  $\Omega_i$.
\end{enumerate}
\edefe

\brem\label{ramindexandtheta} Recall that choice of elements ${\boldsymbol\theta} = \{\theta_i\} \in (Y(T) \otimes \bq)^m$ in the interior of the facets  $\Omega_i$ , identifies the group scheme ${\mathcal G}_{_{\Omega,X}}$ with ${\mathcal G}_{_{{\boldsymbol\theta},X}}$. Starting with a tuple of weights ${\boldsymbol\theta} \in (Y(T) \otimes \bq)^m$, following the proof of  the converse in Theorem \ref{weilbt}, we get positive integers $d_1, d_2, \ldots d_m$ such that $d_i.\theta_i \in Y(T)$. Fix  ${\mathcal R}\subset X$ a finite subset with $|{\mathcal R}| = m$ where the group scheme ${{\mathcal G}_{_{{\boldsymbol\theta},X}}}$ is the local Bruhat-Tits group scheme, with weights $\theta_i$ in the interior of $\Omega_i$.
By choosing the $d_i$ to be the {\em least with this property}, we see that a choice of $\boldsymbol\theta$ entails a choice of ramification indices $d_i$ at  the points of $\mathcal R$.  Then by generalities on ramified covers (see \ref{ramindices}), we can get a covering $p:Y \to X$, ramified over ${\mathcal R}$, with ramification indices $d_i$ and with Galois group $\Gamma$. Note however, that the  local data of $\{d_i\}$ and the ramification points associated to the weights are intrinsic, i.e., depends only on X (see \ref{ramindices}). \erem

\brem The weights $\boldsymbol\theta$ can always be chosen in ${\mathcal A}^m$, where $\mathcal A$ is the Weyl alcove. \erem

\brem\label{classicalwts} The notion of weight defined above  is the precise analogue of the classical weight for a parabolic vector bundle with {\em multiplicity} when the weights are rational(cf. \cite[Definition 1.5, page 211]{ms}). This can be seen by considering  Example \ref{example} which is in fact the original  motivation for parabolic weights. In this context, we refer the reader to Boalch \cite{boalch} where weights come up in a slightly different context. \erem

\subsection{\bf Parabolic line bundles}
Fix a finite subset ${\mathcal R} \subset X$ with $|{\mathcal R}| = m$.
\bdefe\label{parablb} (see \cite[Definition 1.5, page 211]{ms})  A {\em parabolic line bundle} on $(X, {\mathcal R})$ is a pair $(\cl, \{\alpha_1, \ldots,\alpha_m\})$,  where $\cl$ is a line bundle on $X$  together with a $m$--tuple of rational numbers $(\alpha_1, \ldots,\alpha_m)$ with $0 \leq \alpha_i \leq 1$. The {\em parabolic degree} of a parabolic line bundle is defined as
\[
pardeg(\cl) = deg(\cl) + \sum_{i=1}^{m} \alpha_i
\]
\edefe
\brem\label{invdirectimage} Let $p:Y \to X$ be a Galois cover ramified over ${\mathcal R}
 \subset X$ with ramification indices $n_{y_i}, i = 1, \ldots ,m$ at the points $y_i \in Y$ over ${\mathcal R}
$ and let $Gal(Y/X) = \Gamma$. 

Let $L$ be a a $\Gamma$--line bundle on $Y$ of local type $\boldsymbol\tau = \{\boldsymbol\tau_i\}$, where each $\boldsymbol\tau_i$ acts a character $\boldsymbol\tau_i(\zeta) = \zeta^{a_{y_i}}$ with $|a_{y_i}| < n_{y_i}, \forall i$. Then by \cite{pibundles} and \cite{ms}, the {\em invariant direct image} $\cl \simeq p^{\Gamma}_{_*}(L)$ determines a parabolic line bundle on $(X, {\mathcal R}
)$ with parabolic weights $({\frac{a_{y_1}}{n_{y_1}}}, \ldots, {\frac{a_{y_m}}{n_{y_m}}})$ and parabolic degree:
\[
pardeg(p^{\Gamma}_{_*}(L)) = deg(p^{\Gamma}_{_*}(L)) + \sum_{i=1}^{m} {\frac{a_{y_i}}{n_{y_i}}}   
\]
\erem
\brem In fact, all parabolic line bundles on $(X, {\mathcal R})$ can be realized in this manner namely, as {\em invariant direct images}; this is done by constructing a cover ramified over ${\mathcal R}$ with suitable ramification indices. 
\erem

\subsection{\bf Semistability and stability of torsors} Let ${\mathcal G}_{_{\Omega,X}}$ be a Bruhat-Tits group scheme on the curve $X$ as in Definition \ref{globalbt}. Let  $(\eE,\boldsymbol\theta)$ be a parahoric torsor, i.e., the weights $\boldsymbol\theta$ are such that ${\mathcal G}_{_{\Omega,X}} \simeq {{\mathcal G}_{_{{\boldsymbol\theta},X}}}$.     

Let ${P}_K \subset {\mathcal G}_K$ be a maximal parabolic subgroup of the generic fibre ${\mathcal G}_K$ of ${\mathcal G}_{_{\Omega,X}}$. Let $\chi: { P}_K \to {\mathbb G}_{m,K}$ be a dominant character of the parabolic subgroup ${ P}_K$. Then one knows that this defines an ample line bundle $L_{_\chi}$ on ${\mathcal G}_K/{ P}_K$.  We see immediately that $\chi$ defines a line bundle $L_{_\chi}$ on ${\eE}_K({\mathcal G}_K/{ P}_K) \simeq {\eE}_K/{ P}_K$ as well and using a reduction section $s_K$, we therefore get a line bundle $s_K^{*}(L_{_\chi})$ on $X - {\mathcal R}$. 

\bprop\label{better} Let ${\mathcal G}_K$ be the generic fibre of the Bruhat-Tits group scheme ${\mathcal G}_{_{\Omega,X}}$. Let $(\eE,\boldsymbol\theta)$ be a parahoric torsor. Let $s_K$ be a generic reduction of structure group of ${\eE}_K$ to ${ P}_K$. Then the line bundle $s_K^{*}(L_{_\chi})$ on $X - {\mathcal R}$ has a canonical extension $L_{_\chi}^{\boldsymbol\theta}$ to $X$ as a parabolic line bundle. \eprop

\begin{proof} The choice of $\theta \in Y(T) \otimes \bq$ allows us to choose an integer $d$ such that $d.\theta \in Y(T)$. Then we have a ramified cover $p:Y \to X$ with $\Gamma = Gal(Y/X)$ ramified over $x$ with ramification index $d$.  By Theorem \ref{stackygammaversusparahoric}, the parahoric torsor $(\eE,\boldsymbol\theta)$, comes from a $(\Gamma,G)$--principal bundle $E$ of local type $\boldsymbol\tau$ on $Y$; more precisely, $\eE \simeq p^{\Gamma}_{_*}(E \wedge^G F^{op})$ for a fixed $(\Gamma,G)$--bundle $F$. 

The maximal parabolic subgroup ${ P}_K \subset {\mathcal G}_K$ immediately gives a maximal parabolic $Q \subset G$ and the reduction $s_K$ gives in turn a $\Gamma$--equivariant reduction of structure group $t_L$, i.e., a section of  $E_L/Q_L$, where $L$ denotes the quotient field of $B$ the local ring in $Y$ over $x \in X$. By virtue of the projectivity of $Y$, the reduction section $t_L$ extends to a $\Gamma$--equivariant reduction of structure group $t$ as a section of $E/Q$. The dominant character $\chi$ gives a dominant character $\eta$ of $Q$ and the section $t$ gives a $\Gamma$--line bundle $t^{*}(L_{_\eta})$ on $Y$.

 We observe that the line bundle $L_{_\chi}^{\boldsymbol\theta}:= p^{\Gamma}_{_*}(t^{*}(L_{_\eta}))$
gives the required  extension of $s_K^{*}(L_{_\chi})$. By the very definition of the invariant direct image  (see Remark \ref{invdirectimage}), we see that~$L_{_\chi}^{\boldsymbol\theta} = p^{\Gamma}_{_*}(t^{*}(L_{_\eta}))$ gets the natural structure of a {\em parabolic line bundle}. 

Note that the parabolic line bundle extension $L_{_\chi}^{\boldsymbol\theta}$ obtained above depends only on local data coming from $\boldsymbol\theta$ (see Remark \ref {ramindexandtheta}) and hence is intrinsic on $X$, i.e., it does not depend on the choice of $Y$ (see Remark \ref{localdata}).

\end{proof}

We have the following general definition of stability and semistability for $(\Gamma,G)$--bundles following A. Ramanathan \cite[Lemma 2.1]{Ramanathan}.

\bdefe\label{ramstability}( Semistability and stability) Let $G$ be a reductive algebraic group. A   $(\Gamma,G)$--bundle $E$ on $Y$ is called  
$\Gamma$-\textit{semistable} (resp. $\Gamma$-\textit{stable}) if for every maximal parabolic subgroup $P \subset G$ and every $\Gamma$--invariant reduction of structure group $\sigma:Y \to E(G/P)$, and for every dominant character $\chi:P \to {\mathbb G}_m$   we have $\deg\ \sigma^{*}(L_{_{\chi}}) \leq 0.$ (resp $ < 0$).\edefe

\brem\label{equivalentss} (cf. \cite[Definition 1.1]{ram1})~ Equivalently, $E$ is called  
$\Gamma$-\textit{semistable} (resp. $\Gamma$-\textit{stable}) if for every maximal parabolic subgroup $P \subset G$ and every $\Gamma$--invariant reduction of structure group $\sigma:Y \to E(G/P)$, we have $\sigma^{*}(E(\mathfrak{g}/\mathfrak{p})) \geq 0$ (resp $ > 0$), where $\mathfrak{g}$ (resp.$\mathfrak{p}$) is the Lie algebra of $G$ (resp. $P$). Note that $E(\mathfrak{g}/\mathfrak{p})$ can be identified with $E(T_{_{G/P}})$, where $T_{_{G/P}}$ is the relative tangent sheaf to the morphism $E(G/P) \to X$. \erem 

We therefore have the following analogous definition:
\bdefe\label{stability} Let $\mathcal G = {\mathcal G}_{_{\Omega,X}}$. A parahoric torsor $(\eE, {\boldsymbol\theta})$ is called {\em stable (resp. semistable)} if for every maximal parabolic ${\mathcal P}_K \subset {\mathcal G}_K$, for every dominant character $\chi$ as above and for every reduction of structure group $s_K$, we have:
\[
par deg(L_{_\chi}^{\boldsymbol\theta}) < 0 (resp. \leq 0)
\]
\edefe

The following theorem is immediate from the above discussions together with  Definition \ref{ramstability}. 
\bth\label{stableparahoricgamma} The isomorphism
\[
p^{\Gamma}_{_*}:{\sf Bun}_{_Y}^{\boldsymbol\tau}({\Gamma},G) \xrightarrow{\sim} {\sf Bun}_{_X}({{\mathcal G}_{_{{\boldsymbol\theta,X}}}}) 
\]
given by Theorem \ref{stackygammaversusparahoric} identifies the substacks of stable (resp. semistable) parahoric torsors  with stable (resp. semistable) $(\Gamma,G)$--bundles of local type $\boldsymbol\tau$ on the ramified cover $Y$. \eeth

\brem Recall the classical definition of a stable parabolic vector bundle as given in \cite{ms}. Note that the definition in \cite[Definition 1.13]{ms} is the one which arises out of interpreting the $\pi$--stability of the $\pi$--vector for the invariant direct image. By Remark \ref{classicalwts} the notion of parabolic weights defined in \cite{ms}  is the same as the one given here when $G = GL(n)$. Our definition of $\Gamma$--stability for $(\Gamma,G)$--bundles generalizes the one given by A.Ramanathan for $G$--bundles, which generalizes the usual notion of stability of vector bundles.  \erem

\brem\label{parabolicforparahoric} Following \cite[Definition 17]{heinloth} one can define a parabolic subgroup $\mathcal P \subset {\mathcal G}_{_{\theta,X}}$ of the group scheme ${\mathcal G}_{_{\theta,X}}$  as the flat closure of a  parabolic subgroup of the generic fibre ${\mathcal G}_K$ of ${\mathcal G}_{_{\theta,X}}$. 
Let $E$ be a ${\mathcal G}_{_{\theta,X}}$--torsor on $X$. Then as in (\cite[Lemma 23]{heinloth}) one can show that if ${\mathcal P}_K \subset {\mathcal G}_K$ is a parabolic subgroup and $E$ a ${\mathcal G}_{_{\theta,X}}$--torsor on $X$, then any choice of a reduction section $s_K \in E_K({\mathcal G}_K/{\mathcal P}_K)$ defines a parabolic subgroup ${\mathcal P}' \subset {\mathcal G}_{_{\theta,X}}$ together with a reduction $s'$ of $E$ to ${\mathcal P}'$. 

In fact, these results in \cite{heinloth} can be deduced by using the invariant direct image concept. Let $H \subset G$ be a closed subgroup and let $F$ be a fixed $(\Gamma,G)$--bundle of type $\boldsymbol\tau$ with a $\Gamma$--invariant reduction of structure to $H$. Let the induced $(\Gamma,H)$--bundle obtained from this reduction be denoted by $F_{_H}$. We consider the adjoint group scheme ${\mathcal G}_{_{F_{_H}}}$ as defined in Notation~ \ref{bigdelta}. Then, ${\mathcal G}_{_{F_{_H}}} \subset {\mathcal G}_{_F}$ is a closed subgroup scheme. By \cite[Proposition 2, page 192]{blr} and taking $\Gamma$--invariants,  it follows immediately  by \cite[Proposition 3.4]{bass} that $p^{\Gamma}_{_*}({\mathcal G}_{_{F_{_H}}}) \subset p^{\Gamma}_{_*}({\mathcal G}_{_F}) \simeq {\mathcal G}_{_{\theta,X}}$ is a closed smooth $X$--subgroup scheme. In particular, if $P$ is a parabolic subgroup of $G$, the invariant direct image $p^{\Gamma}_{_*}({\mathcal G}_{_{F_{_P}}})$ gives the flat closure $\mathcal P \subset {\mathcal G}_{_{\theta,X}}$. 

\erem

\brem ({\em Harder-Narasimhan reduction}) With the definition of semistability in place, it is routine now to define the Harder-Narasimhan reduction for a ${\mathcal G}_{_{\theta,X}}$--torsor by using the identification of Theorem \ref{stableparahoricgamma}. The existence of a parahoric Harder-Narasimhan reduction follows from the existence of a $\Gamma$--equivariant parabolic Harder-Narasimhan reduction for a $(\Gamma,G)$--bundle together with Remark \ref{parabolicforparahoric}. In other words, the canonical Harder-Narasimhan parabolic subgroup scheme of the parahoric group scheme will be the invariant direct image of the $\Gamma$--invariant Harder-Narasimhan parabolic subgroup $P$ of $G$.  The well-definedness follows since it is the flat closure of the generic Harder-Narasimhan parabolic $P_K \subset {\mathcal G}_K$. The Harder-Narasimhan reduction of structure group for the torsor will be realized as the invariant direct image of the corresponding $(\Gamma,P)$--reduction on $Y$. The uniqueness of the Harder-Narasimhan reduction for $(\Gamma,G)$--bundles shows the uniqueness of the Harder-Narasimhan reduction of a parahoric torsor as well.
\erem

\section{Unitary representations of $\pi$}

\subsection{\bf  Manifold of irreducible unitary representations of $\pi$}
\label{p172S3CI}
Notations in this section are as in the introduction. 

Let $\rho$ be a representation of $\pi$ on a vector space $V$ (over $\br$) such that $d=\dim V$ and let $\rho$ act unitarily.
We now recall the following result from Weil \cite[Page 156]{weil2}, noting that since $\pi$ acts unitarily, it leaves a non-degenerate form on $V$ invariant and therefore in Weil's notation, $i = i' = \dim_{\br} H^0(\pi, \rho)$.
\bprop\label{p172prop7CI}
We have the following equality of dimensions:
\beqa
\dim_{\br}H^1(\pi, \rho) = 2d (g-1) +2 \dim_{\br} H^0(\pi, \rho)+\sum^m_{\nu=1}e_{\nu}
\eeqa
where $e_{\nu}$ is the rank of the endomorphism $(I - \rho (C_{\nu}))$ of $V.$
\eprop

Let $K_G$  be a maximal compact subgroup of $G$  and $Lie(K_G)$ denote {\em the Lie algebra} of $K_G$, which is a real vector space of dimension $d$, where $d = dim(G)$.

As in the introduction, we assume that $X = {\mathbb H}/{\pi}$, with $x \in X$ corresponding to $z \in {\mathbb H}$. Let $\pi_{z}$ be the stabilizer at $z$ (cyclic of order $n_x$) and let $\gamma$ be a generator of $\pi_{z}$ and now let $\rho: \pi \to K_G$ be a unitary representation of $\pi$ (see Definition \ref{localtypeofrep}).

\subsection {\bf Explicit computation when $G$ is simple:} Let $\alpha \in S$, and $\theta_{\alpha}$ be as in \eqref{thetaalpha}; let $\rho_{_{\theta_{\alpha}}}$ be the local representation as in Notation \ref{rhoalpha}. Let $\rho_{_{\theta_{\alpha}}}(\gamma) \in K_G$ be the image of the generator $\gamma$ of $\pi_{z}$. Note that the choice of the simple root $\alpha$ and identification of the representation $\rho$ with $\rho_{_{\theta_{\alpha}}}$  amounts to {\em fixing the local type} of the representation $\rho: \pi \to K_G$, i.e. the {\em conjugacy class} of $\rho(\gamma)$ in $K_G$. 

We denote by Ad $\rho_{_{\theta_{\alpha}}}$, the adjoint transformation on $Lie(K_G)
$, namely if $M \in Lie(K_G)
, \ M\mapsto \rho_{_{\theta_{\alpha}}}(\gamma) M \rho_{_{\theta_{\alpha}}}(\gamma)^{-1}.$ Then we have:
\bprop\label{weilseshu} 
Let ${\text e}({\theta_{\alpha}})$ denote the rank of $(Id - Ad~\rho_{_{\theta_{\alpha}}})$ on $Lie(K_G)
$. Then \small
\beqa\label{explicitexpr}
 e({\theta_{\alpha}}) = dim_{\br}(K_G) - 2 \mu(\alpha) - 2 \nu(\alpha) - {\ell}  
= 2.(dim_{_\bc}(G/P_{_\alpha})) - \mu(\alpha))
\eeqa
where $P_{_\alpha}$ is the maximal parabolic subgroup of $G$ associated to $\alpha$ and 
\beqa\label{mualpha}
\mu(\alpha) = \#\{r \in R^{+} \mid r = c_{\alpha}.\alpha +  \sum_{\beta \neq \alpha} x_{\beta}.\beta \}
\eeqa
\beqa\label{nualpha}
\nu(\alpha) = \#\{r \in R^{-} \mid r ~involves~~simple~~roots \neq \alpha \}
\eeqa
and $\ell = ~ \mid S \mid$.
\eprop

\begin{proof} Make $K_G$ operate on itself by inner conjugation. Then,  rank of (Id-Ad $\rho_{_{\theta_{\alpha}}})$ acting on the Lie algebra $Lie(K_G)
$ equals  {\em the dimension of the orbit through $\rho_{_{\theta_{\alpha}}}(\gamma)$ for the action of $K_G$ on itself by inner conjugation}.

We may assume for the purpose of this computation that $\rho_{_{\theta_{\alpha}}}(\gamma)$ lies in the maximal torus.  We firstly compute the number of roots $r \in R$ so that the corresponding root group $U_r(B)$ is centralized by $\rho_{_{\theta_{\alpha}}}(\gamma)$. Recall from Definition \ref{rhoalpha} that the action of $\rho_{_{\theta_{\alpha}}}(\gamma)$ on $U_r$ is given as follows:
\beqa
\rho_{_{\theta_{\alpha}}}(\gamma). U_r(B) . \rho_{_{\theta_{\alpha}}}(\gamma)^{-1} = U_r(\zeta^{r(\Delta_{_\alpha})} B)
\eeqa
where as seen earlier, $r(\Delta_{_\alpha}) = d. (\theta_{_\alpha}, r)$. Since $\zeta$ is a primitive $d^{th}$--root of unity, we need to compute the $\#\ \{r \in R \mid (\theta_{_\alpha}, r) = \pm 1~~ or~~ 0 \}$. It is easy to see that
\beqa
\{r \in R \mid (\theta_{_\alpha}, r) = \pm 1~~ or~~ 0\} =  \bigcup_{i = 1}^4 A_i(\alpha) 
\eeqa
where for $i = 1,2$, 
\beqa
A_i(\alpha) = \{r \in R^{\pm} \mid r = \pm c_{\alpha}.\alpha +  \sum_{\beta \neq \alpha} \pm x_{\beta}.\beta \}
\eeqa
\beqa
A_3(\alpha) = \{r \in R^{-} \mid r ~involves~~simple~~roots \neq \alpha \}
\eeqa
and
\beqa
A_4(\alpha) = \{r \in R^{+} \mid r ~involves~~simple~~roots \neq \alpha \}
\eeqa
Since the maximal torus centralizes $\rho_{_{\theta_{\alpha}}}(\gamma)$, we see that the dimension of the centralizer of $\rho_{_{\theta_{\alpha}}}(\gamma)$ is
\beqa\label{abovenumber}
\# \{r \in R \mid (\theta_{_\alpha}, r) = \pm 1 ~~or~~ 0\} ~ +~ \mid S \mid
\eeqa
Observe that $\mid A_4 \mid = \mid A_3 \mid$ and $\mid A_1 \mid = \mid A_2 \mid$ . To compute the rank of (Id - Ad $\rho_{_{\theta_{\alpha}}})$, we simply subtract the above number \eqref{abovenumber} from the $dim_{_\br}(K_G)$ to get the first expression for $e(\alpha)$. We see that 
\beqa\label{nualpha1}
\nu(\alpha) = dim_{_\bc}({P}_{_\alpha}/B)
\eeqa
where $P_{_\alpha}$ is the maximal parabolic subgroup of $G$ defined by the simple root $\alpha \in S$. Thus, 
\[ dim_{_\br}(K_G)  - 2.\nu(\alpha) - \ell = dim_{_\bc}(G)  - 2.\nu(\alpha) - \ell = 2. dim_{_\bc}(G/P_{_\alpha}). 
\]
since $2. dim(B) - \ell = dim(G)$.

Hence, $e({\theta_{\alpha}}) = 2.(dim_{_\bc}(G/P_{_\alpha})) - \mu(\alpha))$ and the proposition now follows. \end{proof}

\bcor\label{hyperspecialdim} Let $\alpha \in S$ be such that ${\mathcal P}_{_{\theta_{\alpha}}}(K)^{hs}$ is a  maximal parahoric subgroup in $G(K)$ {\em which is hyperspecial}. Then $e(\theta_{\alpha}) = 0$ and conversely.
\ecor
\begin{proof} By Bruhat-Tits theory, the hyperspecial parahorics are simply the maximal parahorics $\{ {\mathcal P}_{_{\theta_{\alpha}}}(K) \mid \forall \alpha \in S, {\text {with}}~c_{\alpha} = 1 \}$ upto conjugacy by $G(K)$. In these cases,
the number $\mu(\alpha)$ will now be 
\[
\mu(\alpha) =  \#  \{r \in R^{+} \mid r ~{\text {involves}}~ \alpha \}
\]
since the largest possible coefficient for such an $\alpha$ in any positive root is $1$. Hence $\alpha$ is hyperspecial if and only if $\mu(\alpha) = dim(G/P_{\alpha})$ and we are through by the Proposition \ref{weilseshu}. \end{proof}

\subsection {\bf The moduli dimension}~~Let $G$ be semisimple and simply connected.
\bcor  Let $\theta \in {\mathbb E}$ be an arbitrary element in the affine apartment ${\mathbb E}$ and let $\rho_{{_\theta}}$ be the representation defined in Definition \ref{rhoalpha}. Let ${\text e}({\theta})$ denote the rank of $(Id - Ad~\rho_{_{\theta}})$ on $Lie(K_G)
$. Then, 
\beqa\label{etheta}
e(\theta) = dim_{\br}(K_G) - \mid S \mid - 
\# \{r \in R \mid (\theta, r) = \pm 1 ~~or~~ 0\}  
\eeqa
\ecor
\begin{proof} The proof is immediate from the above discussion. Note that when $\theta = \theta_{\alpha}$, the number $e(\theta)$ gets the explicit expression \eqref{explicitexpr}.\end{proof}

Let ${\boldsymbol\tau} = \{{\boldsymbol\tau}_i\}$ be a set of conjugacy classes and let $\boldsymbol\theta_{\boldsymbol\tau} = \{\theta_i\} \in {\mathbb E}^m$ the corresponding set of points of the product of the affine apartments, with $m = |{\mathcal R}
|$.
\bth\label{realdimension} The subset  $R_o \subset R^{\boldsymbol\tau}(\pi, K_G)$ of {\em irreducible representations} is open and non-empty and is further smooth of real dimension equal to
\beqa
(2g -1) dim(K_G) + \sum_{i = 1}^{m} e({\boldsymbol\theta_{\boldsymbol\tau}}).
\eeqa
Let $K_G$ act on $R^{\boldsymbol\tau}(\pi, K_G)$ by inner conjugation. Let ${\overline K_G} = K_G /{centre}$. Then the equivalence classes of irreducible representations corresponds to the quotient space $R_o/{\overline K_G}$; further, there is an open subset $U$ of $R_o$ where ${\overline K_G}$ acts freely; if $U$ is non-empty then the quotient $U/{\overline K_G}$ has the natural structure of a {\em real analytic manifold}  of real dimension 
\beqa\label{modulidimension}
 2. dim_{_\bc}(G) (g - 1) +  \sum_{i = 1}^{m}e({\boldsymbol\theta_{\boldsymbol\tau}})
\eeqa
\eeth

\begin{proof} We follow the arguments in Narasimhan-Seshadri \cite[Proposition 9.2]{ns} or Seshadri \cite[Page 180]{pibundles}. Let $W = \prod W_{_i}$, where $W_{_i}$ is the conjugacy class defined by ${\boldsymbol\tau}_i$. Observe that the group $\pi$ is given by generators and relations as in \eqref{intro5} and the space $R^{\boldsymbol\tau}(\pi, K_G)$ can be identified with the inverse image of identity under the analytic map $\chi: K_G \times \ldots K_G \times W \to K_G$ given by $\chi(a_1, \ldots, a_g, b_1, \ldots, b_g, c_1, \ldots, c_m) = \prod[a_i,b_i].c_1 \ldots c_m$. As in \cite{ns} or \cite{pibundles}, the kernel of the differential of $\chi$ at $\rho$ is given by $Z^1(\pi, Ad~\rho)$. Also the differential is of maximal rank at $\rho$ if and only if $\rho$ is irreducible. Now using Proposition \ref{weilseshu} the theorem follows  as in loc.cit. \end{proof}

\brem It will be shown in Section 7 that the above open subset is non-empty and gets identified with the $\Gamma$--stable bundles whose automorphisms are trivial; furthermore (see Corollary \ref{narasim}), the quotient $R_o/{\overline K_G}$ in fact gets the structure of a {\em complex analytic orbifold} (i.e., with at most finite quotient singularities) of dimension 
\beqa\label{complexmodulidimension}
dim_{_\bc}(R_o/{\overline K_G}) = dim_{_\bc}(G) (g - 1) +  \sum_{i = 1}^{m}{\frac{1}{2}}e({\boldsymbol\theta_{\boldsymbol\tau}})
\eeqa
\erem

\section{The moduli space of parahoric torsors and the main theorem}
The aim of this section is to construct  the moduli space of semistable $(\Gamma,G)$--bundles on $Y$ of local type $\boldsymbol\tau$ (see Definition \ref{ramstability}), or equivalently, by Theorem \ref{stableparahoricgamma}  the moduli space of semistable and stable parahoric torsors. We essentially follow the strategy of Balaji-Seshadri \cite{basa} and Balaji-Biswas-Nagaraj \cite{bbn1}. We briefly outline a proof of  \cite[Theorem 5.8]{bbn1}.

We fix a faithful representation $G \hra GL(n)$ and consider the subscheme of a suitable ``Quot"-scheme parametrizing $\Gamma$--vector bundles on the curve $Y$ which are $\Gamma$--semistable   of local type ${\boldsymbol\tau}$ and we denote this scheme by $Q^{\boldsymbol\tau}_{_{(\Gamma,GL(n))}}$ (see \cite{pibundles} for details where this space is denoted $R^{{\boldsymbol\tau},ss}$). We may equivalently view the points in $Q^{\boldsymbol\tau}_{_{(\Gamma,GL(n))}}$ as $\Gamma$--semistable principal  $(\Gamma,GL(n))$--bundles of local type ${\boldsymbol\tau}$.

 We then define the scheme $Q^{\boldsymbol\tau}_{_{(\Gamma,G)}}$ as the {\em space of $\Gamma$--equivariant reductions of structure group} of the bundles in $Q^{\boldsymbol\tau}_{_{(\Gamma,GL(n))}}$  which consists of those $(\Gamma,G)$--bundle which are of local type ${\boldsymbol\tau}$. It is standard to show that $Q^{\boldsymbol\tau}_{_{(\Gamma,G)}}$  has the {\em local universal property} for families of $\Gamma$--semistable $(\Gamma,G)$--bundles of local type. 
 
 We now use the results in \cite{pibundles} which shows that there is an action of a certain reductive group ${\mathcal H}$ on $Q^{\boldsymbol\tau}_{_{(\Gamma,GL(n))}}$  and the good quotient $M_{_Y}^{\boldsymbol\tau}(\Gamma,n):= Q^{\boldsymbol\tau}_{_{(\Gamma,GL(n))}}//{\mathcal H}$ exists and gives a coarse moduli scheme for the functor of equivalence classes of $\Gamma$--semistable principal $(\Gamma, GL(n))$--bundles on $Y$ of local type $\boldsymbol\tau$. 
 
The map $Q^{\boldsymbol\tau}_{_{(\Gamma,G)}} \to 
Q^{\boldsymbol\tau}_{_{(\Gamma,GL(n))}}$ obtained by taking extension of structure groups via the inclusion $G \hra GL(n)$, is shown to be {\em affine} and the action of ${\mathcal H}$ lifts to $Q^{\boldsymbol\tau}_{_{(\Gamma,G)}}$ to give a good quotient $Q^{\boldsymbol\tau}_{_{(\Gamma,G)}}//{\mathcal H}$ which {\em we denote by} $M_{_Y}^{\boldsymbol\tau}(\Gamma, G)$ (see \cite{basa} and \cite{bbn1}).

When $G$ is semisimple and simply connected, we show in this paper  that the points of the scheme $M_{_Y}^{\boldsymbol\tau}(\Gamma, G)$ parametrize {\em isomorphism classes of $(\Gamma,G)$--bundles of local type $\boldsymbol\tau$ which are unitary (Definition \ref{moreunitarybundle} below)}. Using this we show that $M_{_Y}^{\boldsymbol\tau}(\Gamma, G)$ is normal and projective and compute its dimension.  

\brem We note that the arguments of \cite{bbn1} are not sufficient for showing the last statement (i.e the projectivity and dimension computation) since the local type of the bundles was not fixed in \cite{bbn1}. A key step in the arguments is the connectedness of the moduli space which fails if the local type is not fixed.\erem

\brem Note that strictly speaking, we do not need the group $G$ to be semisimple and simply connected but we need only the reductivity of $G$ to be able to talk of the space $M_{_Y}^{\boldsymbol\tau}(\Gamma, G)$. \erem

\bdefe\label{unitarybundle} A  {\em unitary} $(\pi,G)$--bundle on $\bh$  is defined to be the trivial $G$--bundle $\bh \times G$ on $\bh$ with the  $\pi$--structure given by $\gamma(z,g) = (z, \rho(\gamma).g)$, with $\rho$ an element of $R^{\boldsymbol\tau}(\pi, K_G)$.\edefe  
 
 Let $V$ be a {\em unitary} $(\pi,G)$--bundle defined by $\rho:\pi \to K_{_G}$. Let $r:\bh \to Y$ be as in \eqref{intro1}. Let $E(\rho):= r_{_*}^{\pi_o}(V)$; then $E(\rho)$ is a $(\Gamma,G)$--bundle defined by the {\em twisted action} given by \eqref{intro2}. 
 
 We observe that the {\em local type $\boldsymbol\tau_i$ of the bundle $E(\rho)$ at $y_i$} in the sense of Definition \ref{earlylocaltype} is equivalently given by the {\em conjugacy class of }$\rho(C_i)$ in $G$. Thus if $\boldsymbol\tau = \{\boldsymbol\tau_i \}$, then we have
\beqa\label{intro7}
\mbox{\parbox{3.75in}{
$\rho$ is of type $\boldsymbol\tau = \{\boldsymbol\tau_i \} \Longleftrightarrow E(\rho)$ is of local type $\boldsymbol\tau$}}
\eeqa

\bdefe\label{moreunitarybundle} A $(\Gamma,G)$--bundle $E$ is called {\em unitary} if $E \simeq E(\rho)$ for a homomorphism $\rho:\pi \to K_{_G}$.\edefe

\subsection{\bf  Properness of the moduli of $(\Gamma,G)$--bundles.}

Let $H = G/Z(G)$, the associated adjoint group. 
Let $\mathfrak h = Lie(H)$. Consider the adjoint representation
$\rho:H \to GL(\mathfrak h)$. It  is clear that $\rho$ is  faithful
 representation.

Fix the representation $\rho:H \hra GL(n)$ (where $n = dim~{\mathfrak h}$) and a maximal
compact ${K_H}$ of $H$ such that ${K_H} \hra
U(n)$. Consider the subscheme $M_{_Y}^{\boldsymbol\tau}(\Gamma,n)^{s} \subset M_{_Y}^{\boldsymbol\tau}(\Gamma,n)$  of stable $(\Gamma,GL(n))$--bundles.

\blem\label{goodunitary} Let $\phi :M_{_Y}^{\boldsymbol\tau}(\Gamma,H) \lr M_{_Y}^{\boldsymbol\tau}(\Gamma,n)$ be the morphism induced 
by the representation $\rho$ and the map of Quot schemes. Let $M_{_Y}^{\boldsymbol\tau}(\Gamma,H)^{o}:= {\phi}^{-1}(M_{_Y}^{\boldsymbol\tau}(\Gamma,n)^{s})$ be the inverse image of the  stable points. Then, $M_{_Y}^{\boldsymbol\tau}(\Gamma,H)^{o}$ (when nonempty), is open and
consists of unitary 
$(\Gamma,H)$--bundles which are $\Gamma$--stable as well. \elem

\begin{proof} We claim that a principal
$(\Gamma,H)$ bundle $E$ is unitary if and only if the associated $(\Gamma, GL(\mathfrak h))$--bundle
$E(\mathfrak h)$ is so. If $E$ is unitary obviously so is $E(\mathfrak h)$. 

We now show the converse. Let $A(\mathfrak h)$ denote the stabilizer of the $GL(\mathfrak h)$--action on the tensor space ${\mathfrak h}^{*} \otimes {\mathfrak h}^{*} \otimes {\mathfrak h}$ at the point $[~,~]$, i.e. the Lie bracket. Since we have assumed that $H$ is of adjoint type it implies that $A(\mathfrak h) = Aut(\mathfrak h)$.

Now assume that  $E(\mathfrak h)$ comes from a unitary representation of $\pi$, then we take the Lie bracket morphism $E(\mathfrak h) \otimes E(\mathfrak h) \to E(\mathfrak h)$. Both $E(\mathfrak h) \otimes E(\mathfrak h)$ and $E(\mathfrak h)$  come from unitary representations of $\pi$ and by {\em local constancy} (\cite[Proposition 1.2]{ms}), morphisms of such bundles are induced by morphisms of $\pi$--modules. It now follows that $E(\mathfrak h)$ gets a reduction of structure group to the group $A(\mathfrak h) = Aut(\mathfrak h)$.

Since $H$ is a connected adjoint group, firstly, $Ad(H) = H$ and secondly, it gets identified with the group of inner automorphisms; thus we have a short exact sequence:
\[
1 \to H \to A(\mathfrak h) \to F \to 1
\]
where elements of $F \simeq A(\mathfrak h)/H$ are the {\em outer automorphisms}.
 Again we have a similar exact sequence of compact groups:
\[
1 \to K_H \to K_{A(\mathfrak h)} \to F \to 1
\]
The bundle $E$ is therefore such that $E(A(\mathfrak h))$ is a unitary bundle and comes from a representation $\bar{\chi}: \pi \to K_{A(\mathfrak h)}$. Furthermore, the extended bundle $E(A(\mathfrak h))(F)$ is trivial since it comes with a section (giving $E$). By composing the representation $\bar{\chi}$ with the map $K_{A(\mathfrak h)} \to F$, we see that the triviality of 
$E(A(\mathfrak h))(F)$ forces the composite to be the trivial homomorphism, implying that $\bar{\chi}$ factors via $\chi: \pi \to K_H$ to give the bundle $E$ (cf. Atiyah-Bott \cite[Lemma 10.12]{AB}).

Now using the main theorem of \cite{pibundles} we see that points of
$M_{_Y}^{\boldsymbol\tau}(\Gamma,n)^{s}$, being stable bundles, are all unitary. Hence by 
the claim above the bundles in the inverse image
$\phi^{-1}(M_{_Y}^{\boldsymbol\tau}(\Gamma,n)^{s})$ are also unitary. 

It follows easily from Remark \ref{equivalentss} (cf. \cite[Remark 2.2]{ram1}), that a $(\Gamma,H)$--bundle is $\Gamma$--stable if and only if the associated Lie algebra bundle $E(\mathfrak h)$ has no $\Gamma$--invariant parabolic subalgebra bundles of degree $ \geq 0$. 
It is now easy to see that a $(\Gamma,H)$--bundle is $\Gamma$--stable if the associated Lie algebra bundle is a $\Gamma$--stable vector bundle since $E(\mathfrak h)$ has no  $\Gamma$--subbundles of degree $\geq 0$, and in particular no $\Gamma$--invariant parabolic subalgebra bundles of degree $\geq 0$. \end{proof}

\bprop\label{goodstablebundles} Assume that $H$ is simple of adjoint type. Let $\rho$ be the adjoint representation of $H$. Then the inverse image of 
$M_{_Y}^{\boldsymbol\tau}(\Gamma,n)^{s}$ by the induced morphism $\phi$ 
is nonempty.\eprop

\begin{proof} Recall that the Fuchsian group $\pi$ can be identified with the group generated by $2g + m$ elements $A_i, B_i, C_i$, modulo relations given by \eqref{intro5} and \eqref{intro5.1}.

So to prove that the inverse image $\phi^{-1}({\sf Bun}_{_Y}^{\boldsymbol\tau}(\Gamma,n)^{s})$ is nonempty,
we need to exhibit a representation 
$\chi \, : \, \pi \, \to \, {K_H}$
such that the composition
\beqa
\rho\circ\chi \, : \, \pi \, \to \, U(n) 
\text{~~is {\it irreducible}}.
\eeqa
Choose elements $h_{1},\cdots ,h_{m} \in {K_H}$ so that they are elements
of order $n_{i}$, where $i = 1,\cdots ,m$ (these correspond to
fixing the {\it local type $\boldsymbol \tau$} of our bundles). 

It is a well--known fact that every 
element of a compact connected real semisimple Lie group is a 
commutator. Another well-known fact is  that there exists a dense subgroup
 $\langle \alpha,\beta \rangle $ of ${K_H}$ generated by two general elements $\{\alpha,\beta \}$ (see for example \cite[Lemma 3.1]{Su}).
Recall that the genus $g \geq 2$ and define the representation $\chi:\pi \to {K_H}$ as follows~:

\footnotesize
\beqa
\chi(A_{1})\, =\, \alpha, \, \chi(B_{1})\, =\, \beta, \,
\chi(A_{2})\, =\, \beta, \, \chi(B_{2})\, =\, \alpha, \,
\eeqa  
\beqa \chi(A_{i})\, =\, a_{i}, \, \chi(B_{i})\, =\, b_{i}, \,
~~for~~ i=3,\cdots ,g, \,
 \chi(C_{j})\, =\, h_{j}, \, ~~and~~ \,
j=1,\cdots ,m
\eeqa

\normalsize
It is clear that $\chi$ gives a representation of the group
$\pi$. Since $H$ is {\em simple}, $\rho$ is irreducible, and the image of $\chi$
contains a dense subgroup, the composition $\rho \circ \chi$
gives an irreducible representation of $\pi$ in the unitary
group $U(n)$. Therefore, it gives a {\it stable}
$\Gamma$--linearized vector bundle, which comes as the extension
of structure group of a $H$--bundle. This completes the proof of
the Proposition. \end{proof}

\bcor\label{5.12} There is a
non-empty Zariski open subscheme $M_{_Y}^{\boldsymbol\tau}(\Gamma,H)^{o}$ of $M_{_Y}^{\boldsymbol\tau}(\Gamma,H)$  consisting  of unitary bundles of local type $\boldsymbol \tau$ which are also $\Gamma$--stable.\ecor

\begin{proof} The Corollary follows from Lemma \ref{goodunitary} and Proposition \ref{goodstablebundles}.  For we observe that since $H$ is semisimple of adjoint type, it can be written as a direct product $\prod H_i$ of simple groups of adjoint type. Now a $(\Gamma,H)$--bundle (resp. unitary) is the same as a product of $(\Gamma,H_i)$--bundles (resp. unitary). Likewise by \cite[Proposition 7.1]{ram1}, a $\Gamma$--stable  $(\Gamma,H)$--bundle  is the same as a product of $\Gamma$--stable  $(\Gamma,H_i)$--bundles. For each factor $H_i$, Lemma \ref{goodunitary} and Proposition \ref{goodstablebundles} applies and  the result follows.\end{proof}

We now return to  $G$ which is as before a semisimple, simply connected algebraic group.

\bprop\label{5.120} The subscheme of  $M_{_Y}^{\boldsymbol\tau}(\Gamma,G)$  the  consisting  of stable unitary bundles of local type $\boldsymbol \tau$ is
{\em non-empty} and contains a Zariski open subset.\eprop

\begin{proof} Let $\eta:M_{_Y}^{\boldsymbol\tau}(\Gamma,G) \to M_{_Y}^{\boldsymbol\tau}(\Gamma,H)$ be the morphism induced by the quotient map $G \to H$. Let $M_{_Y}^{\boldsymbol\tau}(\Gamma,H)^{o}$ be as in Corollary \ref{5.12}. We claim that the required Zariski open subset of $M_{_Y}^{\boldsymbol\tau}(\Gamma,G)$ is 
\beqa
M_{_Y}^{\boldsymbol\tau}(\Gamma,G)^{o}:= \eta^{-1}(M_{_Y}^{\boldsymbol\tau}(\Gamma,H)^{o}).\eeqa
Let $E$ be a $(\Gamma,G)$--bundle in $\eta^{-1}(M_{_Y}^{\boldsymbol\tau}(\Gamma,H)^{o})$. By Corollary \ref{5.12}, the $H$--bundle $E(H)$ comes from a unitary representation $\rho: \pi \to K_H$.

Recall that, by the structure of $\pi$ described above, there is a central extension
\beqa\label{centralextension}
1 \to Z_{\tilde{\pi}} \to \tilde{\pi} \to \pi \to 1
\eeqa
where $\tilde{\pi}$ is generated by $A_1, \ldots, A_g, B_1, \ldots, B_g, C_1, \ldots, C_m$ together with a central element $J$ satisfying the extra relation
\beqa\label{centralextension1}
[A_{1},B_{1}]\cdots [A_{g},B_{g}]\cdot C_{1}
\cdots C_{m} = J.
\eeqa

It is easy (as in \cite{ns}), by adding an extra lasso around a dummy point (other than the parabolic points) to choose a lift of $\rho$ to a representation $\tilde{\rho}: \tilde{\pi} \to K_G$ so that the associated $(\Gamma,G)$--bundle $E(\tilde{\rho})$ also maps to $E(H)$. Thus, both $E$ and $E(\tilde{\rho})$ give $E(H)$ under the quotient map $G \to H$. Therefore, by twisting by a central character of $\tilde{\pi}$, we get a representation $\tilde{\pi} \to K_G$  which gives the $(\Gamma,G)$--bundle $E$ (cf. \cite[Page 148]{ram1}). 

We observe that this representation $\tilde{\pi} \to K_G$ in fact descends to a representation $\pi \to K_G$. This follows from the fact that the local type of $E$ at the dummy point is trivial.

From this we can now conclude that all bundles in $M_{_Y}^{\boldsymbol\tau}(\Gamma,G)^{o}$ are unitary (cf. \cite[Lemma 10.12]{AB}). Furthermore, since $G \to H$ is surjective,  it is not hard to see that a $(\Gamma,G)$--bundle is $\Gamma$--stable if and only if the associated $(\Gamma,H)$--bundle is so (cf. \cite[Proposition 7.1]{ram1}). It follows that all points of $M_{_Y}^{\boldsymbol\tau}(\Gamma,G)^{o}$ are also $\Gamma$--stable $(\Gamma,G)$--bundles, completing the proof of the proposition. \end{proof}

\noindent
We now have  a canonical continuous map
\beqa\label{intro8}
\psi:R^{\boldsymbol\tau}(\pi,K_G) \to M_{_Y}^{\boldsymbol \tau}(\Gamma,G)
\eeqa
which sends $\rho$ to the class of~$E(\rho)$. This map is obtained following \cite[page 334]{unitary}. First we consider the  space $R^{\boldsymbol\tau}(\pi,G)$ of all homomorphisms $\pi\to G$ of local type ${\boldsymbol\tau}$. Let $R^{\boldsymbol\tau}(\pi,G)^{ss}$ be the subset of $R^{\boldsymbol\tau}(\pi,G)$ consisting of points $\rho$ such that $E(\rho)$ is $\Gamma$--semistable. One can easily construct an analytic family of $(\pi,G)$--bundles on $\bh \times R^{\boldsymbol\tau}(\pi,G)$; the subgroup $\pi_o$ acts freely on $\bh$ and this family is easily seen to come down to an analytic family of $(\Gamma,G)$--bundles on $Y$ parametrized by $R^{\boldsymbol\tau}(\pi,G)^{ss}$. 

Further, since a $(\Gamma,G)$--bundle is $\Gamma$--semistable if and only if the associated $\Gamma$--vector bundle is so (see for example proof of \cite[Proposition 3.2]{bbn1}), it follows that the subset $R^{\boldsymbol\tau}(\pi,G)^{ss}$ is non-empty and open in 
$R^{\boldsymbol\tau}(\pi,G)$ and contains the space $R^{\boldsymbol\tau}(\pi,K_G)$ of all unitary representations.

By the local universal property of $Q^{\boldsymbol\tau}_{_{(\Gamma,G)}}$, given a $\rho \in  R^{\boldsymbol\tau}(\pi,G)^{ss}$, we get an analytic neighbourhood $U$ of $\rho$ together with an analytic map $U \to (Q^{\boldsymbol\tau}_{_{(\Gamma,G)}})^{ss}$. These maps glue to give an analytic morphism $\psi:R^{\boldsymbol\tau}(\pi,G)^{ss} \to M_{_Y}^{\boldsymbol \tau}(\Gamma,G)$.  Restricting this map to   
$R^{\boldsymbol\tau}(\pi,K_G)$ gives the continuous map $\psi$. The image of $\psi$ consists of $(\Gamma,G)$--bundles which are {\em unitary}.

The following irreducibility result is an immediate consequence of Theorem \ref{stackygammaversusparahoric},  \cite[Theorem 2]{heinloth} and \cite[Proposition 1]{heinloth}.
\bprop\label{stackirred} The moduli stack ${\sf Bun}_{_Y}^{\boldsymbol\tau}(\Gamma, G)$ of $(\Gamma,G)$--bundles on $Y$ of local type $\boldsymbol\tau$ is irreducible and smooth when the group $G$ is semisimple and simply connected. \eprop

\brem We now indicate  a different proof of the connectedness from the picture of Hecke correspondences shown in \eqref{smallhecke}.
By Drinfeld-Simpson\cite{ds}, the moduli stack ${\sf Bun}_{_X}(G)$ is irreducible because $G$ is semisimple and simply connected. Further, the morphism ${\sf Bun}({\mathcal G}_{_{{\mathcal I},X}}) \to {\sf Bun}_{_X}(G)$ is surjective and has fibre $G/B$, $B$ being the Borel subgroup. Hence, ${\sf Bun}({\mathcal G}_{_{{\mathcal I},X}})$ is connected. Now observe that the map ${\sf Bun}({\mathcal G}_{_{{\mathcal I},X}}) \to {\sf Bun}({\mathcal G}_{_{\Omega,X}})$ given by \eqref{smallhecke} is also surjective since it comes from the inclusion ${\mathcal I} \subset {\mathcal P}_{_\Omega}(K)$. Hence ${\sf Bun}({\mathcal G}_{_{\Omega,X}})$ is connected. The irreducibility follows from the formal smoothness of the functor of torsors (see \cite[Proposition 1]{heinloth}; the obstruction to smoothness vanishes since we work on curves.

Since we work over char $0$, the connectedness of the moduli space of $(\Gamma,G)$--bundles of local type $\boldsymbol\tau$  could  also be carried out following 
\cite{ram1} and \cite[Proposition 4.2]{AB}.
\erem

We have a  morphism $f :{\sf Bun}_{_Y}^{\boldsymbol\tau}(\Gamma,G)^{ss} \to M_{_Y}^{\boldsymbol\tau}(\Gamma, G)$ namely, the canonical 
quotient map obtained by the categorical quotient property of the moduli space $M_{_Y}^{\boldsymbol\tau}(\Gamma, G)$. The map $f$ is surjective on points; therefore by Proposition \ref{stackirred}, this implies that $M_{_Y}^{\boldsymbol\tau}(\Gamma, G)$ is {\it irreducible}.

\bth\label{compact} The map $\psi:R^{\boldsymbol\tau}(\pi, K_G) \to M_{_Y}^{\boldsymbol\tau}(\Gamma, G)$ obtained in \eqref{intro8} is surjective and hence $M_{_Y}^{\boldsymbol\tau}(\Gamma, G)$ is
compact. Further, the variety $M_{_Y}^{\boldsymbol\tau}(\Gamma, G)$ gets a structure of a normal
projective variety. \eeth

\begin{proof} By the Proposition \ref{5.120}, the  subset $M_{_Y}^{\boldsymbol\tau}(\Gamma, G)^{o}$ is nonempty
and consists entirely of unitary bundles. Thus it is a subset of the image $\psi(R^{\boldsymbol\tau}(\pi, K_G))$ in $M_{_Y}^{\boldsymbol\tau}(\Gamma, G)$, i.e., the image
${\psi}(R^{\boldsymbol\tau}(\pi, K_G))$ contains a
Zariski open subset of $M_{_Y}^{\boldsymbol\tau}(\Gamma, G)$. Since $R^{\boldsymbol\tau}(\pi, K_G)$ is
compact the image ${\psi}(R^{\boldsymbol\tau}(\pi, K_G))$  is therefore the whole of
$M_{_Y}^{\boldsymbol\tau}(\Gamma, G)$, 
because  $M_{_Y}^{\boldsymbol\tau}(\Gamma, G)$ is {\it irreducible}. 

This proves that $M_{_Y}^{\boldsymbol\tau}(\Gamma, G)$ is topologically compact
and hence by GAGA it is a projective variety. The normality follows from the smoothness of the stack ${\sf Bun}_{_Y}^{\boldsymbol\tau}(\Gamma,G)^{ss}$, again by Proposition \ref{stackirred}.\end{proof}

\bcor\label{narasim}
\begin{enumerate} 
\item Let $g(X) \geq 2$. Then the map $\psi :R^{\boldsymbol\tau}(\pi, K_G) \to M_{_Y}^{\boldsymbol\tau}(\Gamma, G)$ defined above descends to a map 
\beqa\label{narases}
\psi^{*}:R^{\boldsymbol\tau}(\pi, K_G)/{\overline K_G} \to M_{_Y}^{\boldsymbol\tau}(\Gamma, G)
\eeqa
which gives a homeomorphism of topological spaces. Further, the subset $R_o/{\overline K_G}$ of equivalence classes of irreducible unitary representations maps bijectively onto the subset of stable $(\Gamma,G)$--bundles. 

\item Let $g(X) < 2$. When $X = \bp^1$ and $|\mathcal R| \geq 3$ or when $X$ is an elliptic curve and $\mathcal R \neq \emptyset$, the map $\psi^*$~ in~\eqref{narases} is a homeomorphism provided there exists an irreducible representation $\rho: \pi_1(X - \mathcal R) \to K_G$ with preassigned conjugacy classes  of images of lassos around the points of $\mathcal R$.
\end{enumerate}
\ecor
\begin{proof} The surjectivity of the map $\psi^{*}:R^{\boldsymbol\tau}(\pi, K_G)/{\overline K_G} \to M_{_Y}^{\boldsymbol\tau}(\Gamma, G)$ follows from surjectivity statement in Theorem \ref{compact}. 

For the injectivity of $\psi^{*}$, suppose that $\psi(\rho_1)= \psi(\rho_2)$ i.e, we have an isomorphism $E_{_{\rho_1}} \simeq E_{_{\rho_2}}$ of the unitary  bundles  defined by the $\rho_i$. Now we follow Ramanathan \cite[Proposition 6.2]{ram1} and work with our $\Gamma$ instead of $\pi_1(X - x_o)$. The proof simply goes through and this implies that the $\rho_i$ are in the same orbit of ${\overline K_G}$. One could also argue as in Lemma \ref{goodunitary} to get the  injectivity statement. 

Since $R^{\boldsymbol\tau}(\pi, K_G)/{\overline K_G}$ is compact and $M_{_Y}^{\boldsymbol\tau}(\Gamma, G)$ is Hausdorff in the usual topology, the map $\psi^{*}$ is a homeomorphism. The fact that irreducible representations give stable bundles and vice versa follows exactly as in \cite{ram1}. The second part when the genus is $0$ or $1$ follows from Remark \ref{smallgenus} and Theorem \ref{realdimension}. 
\end{proof}

Let ${\mathcal G}_{_{\Omega,X}}$ be a parahoric Bruhat-Tits group scheme associated to a collection of facets $\Omega = \{\Omega_i\}$.  Choose $\boldsymbol\tau = \{\tau_i\}$ and $\boldsymbol\theta_{\boldsymbol\tau} \in (Y(T) \otimes \bq)^m$, so that ${\mathcal G}_{_{\Omega,X}} \simeq {{\mathcal G}_{_{{\boldsymbol\theta_{\boldsymbol\tau}},X}}}$. Recall that Theorem \ref{stableparahoricgamma} identifies stable (resp. semistable) families of parahoric ${\mathcal G}_{_{\Omega,X}}$--torsors with stable (resp. semistable) $(\Gamma,G)$--bundles of local type $\boldsymbol\tau$ on the ramified cover $Y$.

\bdefe Say two parahoric ${\mathcal G}_{_{\Omega,X}}$--torsors $(E,\boldsymbol\theta)$ and $(F,\boldsymbol\theta)$
on $X$ are $S$--equivalent if the corresponding $(\Gamma,G)$--bundles on $Y$ are $S$--equivalent.\edefe

\brem Recall that notion of $S$--equivalence of principal bundles in \cite{Ramanathan}. It is routine to extend this notion to $(\Gamma,G)$--bundles as well (see \cite{bbn1} and \cite{tw}) . The notions of {\em admissible reduction of structure group} is made with the additional $\Gamma$--equivariance property in \cite{bbn1} and \cite{bbn2}. This gives the analogous definitions of $(\Gamma,G)$--polystable bundles and $\Gamma$--associated graded of a $(\Gamma,G)$--semistable bundle (\cite{tw}). We omit the details.\erem

 Let 
\beqa
M({{\mathcal G}_{_{{\boldsymbol\theta},X}}}) := \left\{
\begin{array}{l}
\mbox{the set of $S$--equivalence classes of } \\ 
\mbox{semistable parahoric ${\mathcal G}_{_{\Omega,X}}$--torsors on $X$}\\ 
\end{array} \right\}\eeqa
and let $M({{\mathcal G}_{_{{\boldsymbol\theta},X}}})^{s} \subset M({{\mathcal G}_{_{{\boldsymbol\theta},X}}})$ denote the subset of {\em stable} torsors.

By definition we have the following set-theoretic identification: 
\beqa
M_{_Y}^{\boldsymbol\tau}(\Gamma, G) \simeq M({{\mathcal G}_{_{{\boldsymbol\theta_{\boldsymbol\tau}},X}}})
\eeqa
and by transport of structure we get  the structure of a variety on $ M({{\mathcal G}_{_{{\boldsymbol\theta},X}}})$. We summarize this discussion in the following theorem which is immediate from  Theorem \ref{compact}: 

\bth\label{maintheoremI,II}  The set $M({{\mathcal G}_{_{{\boldsymbol\theta_{\boldsymbol\tau}},X}}})$ gets a natural structure of an irreducible normal projective variety with $M({{\mathcal G}_{_{{\boldsymbol\theta_{\boldsymbol\tau}},X}}})^{s}$ as an open subset. It gives a {\em coarse moduli space} for the  substack ${\sf Bun}({{\mathcal G}_{_{{\boldsymbol\theta_{\boldsymbol\tau}},X}}})^{ss}$ of semistable torsors in ${\sf Bun}({{\mathcal G}_{_{{\boldsymbol\theta_{\boldsymbol\tau}},X}}})$.  Furthermore, we have a homeomorphism 
\beqa
\phi^{*}:R^{\boldsymbol\tau}(\pi, K_G)/{\overline K_G} \to M({{\mathcal G}_{_{{\boldsymbol\theta_{\boldsymbol\tau}},X}}})
\eeqa
which identifies $R_o/{\overline K_G}$ with $M({{\mathcal G}_{_{{\boldsymbol\theta_{\boldsymbol\tau}},X}}})^{s}$.
\eeth

The next corollary follows  from Theorem \ref{realdimension} and the Theorem \ref{maintheoremI,II}.
\bcor\label{dimensionofmodulispace} Let ${\boldsymbol\theta_{\boldsymbol\tau}} = \{{\theta_i}\} \in {\mathbb E}^m$ be the corresponding point in the product of the affine apartment. Then the dimension of the moduli space  $M({{\mathcal G}_{_{{\boldsymbol\theta_{\boldsymbol\tau}},X}}})$ is given by \beqa
dim_{_{\bc}}(G)(g -1) + \sum_{i = 1}^{m}{\frac{1}{2}}e({\boldsymbol\theta_{\boldsymbol\tau}})
\eeqa
\ecor

\bsem\label{reductivecase}{\it Extension to the case when the structure group is reductive.} We  indicate briefly how to extend the construction of the moduli space of $(\Gamma,H)$--bundles to the case when the structure group $H$ is a connected reductive algebraic group and identify it with the space of homomorphisms from $\pi$ to $K_H$. However, the corresponding relationship with parahoric torsors  for reductive $G$ needs a closer analysis of Bruhat-Tits theory for reductive groups.

 Let $S= [H,H]$ be the derived group, i.e. {\em the maximal connected semisimple subgroup} of $H$. Let $Z_0$ be the connected component of the centre of $H$ (which is a torus) and one know that $S$ and $Z_0$ together generate $H$. Let $G= Z_0 \times S$. Then in fact, $G \to H$ is a finite covering map. It is easy to see (following \cite[page 145]{ram1}) that $(\Gamma,G)$--bundles gives rise to $(\Gamma,H)$--bundles and the stability and semistability of the associated $(\Gamma,H)$--bundles follows immediately from that of the $(\Gamma,G)$--bundles. 

The problem of handling the reductive group $G$ reduces to the problem of handling the semisimple group $H$ but which is {\em not simply connected}. Let $\tilde H$ be the semisimple, simply connected algebraic group which is the covering group of $H$. 

We are in the situation of Proposition \ref{5.120}. Recall the central extension \eqref{centralextension}. By adding a dummy point other than the parabolic point, the theory of $(\pi,H)$--bundles is recovered from that of $(\tilde\pi,\tilde H)$--bundles. Notice that a homomorphism $\pi \to K_H$ has as many liftings $\tilde\pi \to K_{\tilde H}$ as the order of the centre of $\tilde H$. It follows quite easily, following arguments as in Lemma \ref{goodunitary}, that the number of connected components of the moduli space in the non-simply connected case is given by the order of the centre of $\tilde H$. In fact, $Hom(\tilde\pi, K_{\tilde H})$ is a union of spaces labelled by elements of the centre of $\tilde H$. Let $Z_0 = Ker(\tilde H \to H)$. Then there is an action of $H^1(X,Z_0)$ on a specific labelled subset of $Hom(\tilde\pi, K_{\tilde H})$.
A component of the moduli space of representations into $K_H$ can be obtained as a quotient of each of these by the action of $H^1(X,Z_0)$.   Details of these ideas are again found in \cite[page 148]{ram1} and follow the ideas of Narasimhan and Seshadri \cite{ns}, where the data over a dummy point is called a {\em special parabolic structure}. 

\esem

\bsem{\it Hecke Correspondences.} Recall that for the case of linear groups one has the classical Hecke correspondences due to Narasimhan and Ramanan \cite{hecke}. In what follows, we consider parahoric subgroups ${\mathcal P}_{_\Omega}(K)$ of $G(K)$ which contain a fixed Iwahori subgroup ${\mathcal I}$ (see \ref{stdparahorics} for notation). Using \eqref{inclusions}, we get ${\mathcal I} \subset {\mathcal P}^{st}_{_{\alpha}}(K) \subset {\mathcal P}_{_{\theta_{\alpha}}}(K) \cap {\mathcal P}_{_{0}}(K)$. These maps of parahoric groups induce morphisms of the corresponding parahoric Bruhat-Tits group schemes, ${\mathcal G}_{_{{\mathcal I}}} \to {\mathcal G}^{st}_{_{{\alpha}}}$ and ${\mathcal G}_{_{{\mathcal I}}} \to {\mathcal G}_{_{\theta_{\alpha}}}$ and  morphisms at the level of stacks and we obtain the following generalized Hecke correspondences.   The dimension formulae  (see Corollary~\ref{dimensionofmodulispace}) get reflected accurately in the picture.
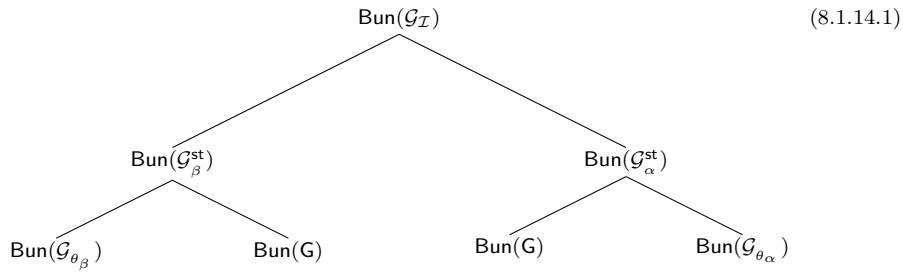
\esem 
\footnotesize
\beqa\label{hecke}
\Tree [.$\sf Bun({\mathcal G}_{\mathcal I})$ [.$\sf Bun({\mathcal G}^{st}_{{_\beta}})$ [.$\sf Bun({\mathcal G}_{_{\theta_{\beta}}})$ ] [.$\sf Bun(G)$ ] ].$\sf Bun({\mathcal G}^{st}_{{_\beta}})$ [.$\sf Bun({\mathcal G}^{st}_{{_\alpha}})$ [.$\sf Bun(G)$ ][.$\sf Bun({\mathcal G}_{_{\theta_{\alpha}}})$ ] ].$\sf Bun({\mathcal G}^{st}_{{_\alpha}})$ ] \eeqa
\\
\normalsize
For instance, we have the following picture of a Hecke correspondence induced by the morphisms ${\mathcal G}_{\mathcal I} \to {\mathcal G}_{_{\Omega}}$ and ${\mathcal G}_{\mathcal I} \to {\mathcal G}_{_{0}} (= G \times X)$:

\footnotesize
\beqa\label{smallhecke}
\xymatrix{
& \sf Bun({\mathcal G}_{\mathcal I}) \ar[dl] \ar[dr]^{G/B} \\
\sf Bun({\mathcal G}_{_{\Omega}}) &   & \sf Bun(G)  
}
\eeqa
\normalsize
\brem It would be interesting to express these relations as morphisms between moduli spaces $M({\mathcal G}_{_{\Omega,X}})$; even the existence of suitable morphisms between the moduli spaces would involve choice of polarization (in the sense of GIT) which would be needed for an algebro-geometric construction of the moduli spaces of parahoric torsors.  \erem


\end{document}